\documentclass[10pt]{amsart}
\usepackage{epsfig}
\usepackage{amsmath, amssymb, euscript,  enumerate}
\usepackage{bbm}
\usepackage{color}
\usepackage{verbatim}
\newcommand{\R}{{\mathbb R}}

\newcommand{\supp}{\text{\rm supp}}

\newcommand{\ap}{\alpha}             
\newcommand{\bt}{\beta}
\newcommand{\gm}{\gamma}             
\newcommand{\dt}{\delta}             
\newcommand{\vep}{\varepsilon}

\newcommand{\ld}{\lambda}            \newcommand{\Ld}{\Lambda}
\newcommand{\sm}{\sigma}             
\newcommand{\vp}{\varphi}
\newcommand{\om}{\omega}             \newcommand{\Om}{\Omega}
\newcommand{\vr}{\varrho}            \newcommand{\iy}{\infty}
            
\newcommand{\f}{\frac}             \newcommand{\el}{\ell}

\newcommand{\fF}{{\mathfrak F}}

\newcommand{\fL}{{\mathfrak L}}

\newcommand{\BI}{{\mathbb I}}

\newcommand{\BN}{{\mathbb N}}

\newcommand{\BR}{{\mathbb R}}

\newcommand{\cB}{{\mathcal B}}

\newcommand{\cF}{{\mathcal F}}

\newcommand{\cK}{{\mathcal K}}

\newcommand{\la}{\langle}          \newcommand{\ra}{\rangle}
\newcommand{\s}{\setminus}         \newcommand{\ep}{\epsilon}
\newcommand{\n}{\nabla}            \newcommand{\e}{\eta}
\newcommand{\pa}{\partial}        
       
    \newcommand{\ds}{\displaystyle}

\newcommand{\bI}{\bold I}        
\newcommand{\bM}{\bold M}        
   
 \newcommand{\pf }{\noindent{\it Proof. }}
\newcommand{\rk }{\noindent{\it Remark. }}

  \newcommand{\pv }{\text{\rm p.v.}}
  \newcommand{\dd }{\text{\rm d}}

\newcommand{\rB }{{\text{\rm B}}}   \newcommand{\rC }{{\text{\rm C}}}
\newcommand{\rL }{{\text{\rm L}}}  
\newcommand{\rP }{{\text{\rm P}}}  \newcommand{\rN }{{\text{\rm N}}}
   \newcommand{\rI}{{\text{\rm I}}}

\newcommand{\btau}{\boldsymbol\tau}

\newtheorem{thm}[subsection]{Theorem}
\newtheorem{lemma}[subsection]{Lemma}
\newtheorem{cor}[subsection]{Corollary}

\newtheorem{definition}[subsection]{Definition}

\numberwithin{equation}{section}

\title[The Evans-Krylov theorem]{The Evans-Krylov theorem for nonlocal \\parabolic fully nonlinear equations}
\author{ Yong-Cheol Kim and Ki-Ahm Lee }
\begin{document}
\begin{abstract} In this paper, we prove the Evans-Krylov theorem
for nonlocal parabolic fully nonlinear equations.
\end{abstract}

\thanks {\bf 10 October, 2015.}

\address{$\bullet$ Yong-Cheol Kim : Department of Mathematics Education, Korea University, Seoul 136-701,
Republic of Korea }

\email{ychkim@korea.ac.kr}

\address{$\bullet$ Ki-Ahm Lee : Department of Mathematics, Seoul National University, Seoul 151-747,
Republic of Korea $\&$ School of Mathematics, Korea Institute for Advanced Study, Seoul 130-722, Republic of Korea}

\email{kiahm@math.snu.ac.kr}

\maketitle

\tableofcontents

\section{Introduction}

Evans and Krylov proved independently an interior regularity
for elliptic partial differential equations which states that any
solution $u\in C^2(B_1)$ of a uniformly elliptic and fully nonlinear
concave equation $F(D^2 u)=0$ in the unit ball $B_1\subset\BR^n$
satisfies an interior estimate $\|u\|_{C^{2,\ap}(B_{1/2})}\le
C\,\|u\|_{C^{1,1}(B_1)}$ with some universal constants $C>0$ and
$\ap\in(0,1)$, so-called the {\it Evans-Krylov} theorem (see
\cite{Ev}, \cite{Kr} and \cite{CS2}). Recently, Caffarelli and
Silvestre \cite{CS1} proved a nonlocal elliptic version of the
Evans-Krylov theorem which describes that any viscosity solution
$u\in L^{\iy}(\BR^n)$ of concave homogeneous equation on
$B_1\subset\BR^n$ formulated by elliptic integro-differential
operators of order $\sm\in(0,2)$ satisfies an estimate
$\|u\|_{C^{\sm+\ap}(B_{1/2})}\le C\,\|u\|_{L^{\iy}(\BR^n)}$ with
universal constants $C>0$ and $\ap\in(0,1)$. This nonlocal
result makes it possible to recover the Evans-Krylov theorem as
$\sm\to 2^-$. In this paper, we prove a parabolic version of the
nonlocal elliptic result of Caffarelli and Silvestre.

We consider the linear {\it parabolic integro-differential
operators} given by
\begin{equation}
\rL u(x,t)-\partial_t
u(x,t)=\pv\int_{\BR^n}\mu_t(u,x,y)K(y)\,dy-\partial_t u(x,t)
\end{equation}
for $\mu_t(u,x,y)=u(x+y,t)+u(x-y,t)-2u(x,t)$. Here we write
$\mu(u,x,y)=u(x+y)+u(x-y)-2u(x)$ if $u$ is independent of $t$. We
refer the detailed definitions of notations to \cite{CS1, KL1,
KL2,KL3}. Then we see that $\rL u(x,t)$ is well-defined provided
that $u\in\rC_x^{1,1}(x,t)\cap\rB(\BR^n_T)$ where $\rB(\BR^n_T)$
denotes {\it the family of all real-valued bounded functions defined
on $\BR^n_T:=\BR^n\times (-T,0]$} and $\rC_x^{1,1}(x,t)$ means
$\rC^{1,1}$-function in $x$-variable at a given point $(x,t)$.
Moreover, $\rL u(x,t)$ is well-defined even for
$u\in\rC_x^{1,1}(x,t)\cap L^{\iy}_T(L^1_{\om})$ (see \cite{KL4}).

We say that the operator $\rL$ belongs to $\fL_0=\fL_0(\sm)$ if its
corresponding kernel $K\in\cK_0=\cK_0(\sm)$ satisfies the uniform
ellipticity assumption:
         \begin{equation}
         (2-\sigma)\frac{\lambda}{|y|^{n+\sigma}}\leq K(y)\leq
         (2-\sigma)\frac{\Lambda}{|y|^{n+\sigma}},\,\,0<\sm<2.
         \end{equation}
If $K(y)=c_{n,\sm}|y|^{-n-\sm}$ where $c_{n,\sm}>0$ is the
normalization constant comparable to $\sm(2-\sm)$ given by
$$c_{n,\sm}=\biggl(\int_{\BR^n}\f{1-\cos(y_1)}{|y|^{n+\sm}}\,dy\biggr)^{-1},$$
then the corresponding operator is
$\rL=-(-\Delta)^{\sm/2}$. Also we say the operator $\rL\in\fL_0$
belongs to $\fL_1=\fL_1(\sm)$ if its corresponding kernel
$K\in\cK_1=\cK_1(\sm)$ satisfies $K\in\rC^1$ away from the origin
and satisfies
         \begin{equation}
        |\n K(y)|\leq \frac{C}{|y|^{n+1+\sigma}}.
         \end{equation}
Finally we say that the operator $\rL\in\fL_1$ belongs to
$\fL_2=\fL_2(\sm)$ if its corresponding kernel
$K\in\cK_2=\cK_2(\sm)$ satisfies $K\in\rC^2$ away from the origin
and satisfies
         \begin{equation}
        |D^2 K(y)|\leq \frac{C}{|y|^{n+2+\sigma}}.
         \end{equation}
The maximal operators are defined by
         \begin{equation*}
         \begin{split}
         \bM_0^+u(x,t)&=\sup_{\rL\in\fL_0}\rL u(x,t)=(2-\sigma)\int_{\R^n}
         \frac{\Lambda\mu^+_t(u,x,y)-\lambda\mu^-_t(u,x,y)}{|y|^{n+\sigma}}dy,\\
         \bM_1^+u(x,t)&=\sup_{\rL\in\fL_1}\rL u(x,t)\,\,\text{ and
         }\,\,\bM_2^+u(x,t)=\sup_{\rL\in\fL_2}\rL u(x,t).
         \end{split}
         \end{equation*}

We shall consider nonlinear integro-differential operators, which
originates from stochastic control theory with jump processes
related with
\begin{equation*}
\bI u(x,t)=\inf_{\bt\in\cB}\rL_{\bt}u(x,t),
\end{equation*} where $\rL_{\bt}
u(x,t)=\pv\int_{\BR^n}\mu_t(u,x,y)K_{\bt}(y)\,dy$ (see \cite{AK,
CS1, KL1, KL2, MP, MR} for the elliptic case and \cite{KL3, KL4} for
the parabolic case). In this paper, we are mainly interested in the
nonlocal parabolic concave equations
\begin{equation}\bI u(x,t)-\pa_t u(x,t)=0\,\,\text{ in
$Q_1$. }\end{equation}

\noindent{[\bf Notations and Definitions]} Let $\sm\in(0,2)$ and
$r>0$.
\begin{itemize}

\item Denote by $Q_r=B_r\times I_r^{\sm}$ and $Q_r(x,t)=Q_r+(x,t)$
for $(x,t)\in\BR^n_T$, where $B_r(x)$ is the open ball with center
$x\in\BR^n$ and radius $r>0$, $B_r=B_r(0)$ and
$I_r^{\sm}=(-r^{\sm},0]$.

\item For a bounded domain $\Om\subset\BR^n$ and $\tau\in(0,T)$, we
denote the parabolic boundary of $\Om_{\tau}=\Om\times(-\tau,0]$ by
$\pa_p\Om_{\tau}:=\pa_x\Om_{\tau}\cup\pa_b\Om_{\tau}:=\pa\Om\times(-\tau,0]\cup\Om\times\{-\tau\}$.

\item The parabolic distance $\dd$ between $X=(x,t)$ and
$Y=(y,s)$ is defined by
\begin{equation}
\dd(X,Y)=
\begin{cases}
(|x-y|^{\sigma}+|t-s|)^{1/\sigma}, &t\le s,\\
\iy, &t>s.
\end{cases}
\end{equation}
For $X_0=(x_0,t_0)\in\BR^n_T$, we set
$\rB^{\dd}_r(x_0,t_0)=\{(x,t)\in\BR^n_T:\dd\bigl(X,X_0)\bigr)<r\}$.

\item We denote by $\om_{\sm}(y)=1/(1+|y|^{n+\sm})$ and $\om:=\om_{\sm_0}$
for some $\sm_0\in(1,2)$ very close to $1$, and also we denote by
$\om(B_{r/2})=\int_{B_{r/2}}\om(y)\,dy$. Let $\fF(\BR^n_T)$ be the
family of all real-valued measurable functions defined on
$\BR^n_T:=\BR^n\times(-T,0]$. For $u\in\fF(\BR^n_T)$ and
$t\in(-T,0]$, we define the weighted norm
$\|u(\cdot,t)\|_{L^1_{\om}}$ by
$$\|u(\cdot,t)\|_{L^1_{\om}}=\int_{\BR^n}|u(x,t)|\om(x)\,dx.$$
Consider the function space $L^{\iy}_T(L^1_{\om})$ of all continuous
$L^1_{\om}$-valued functions $u\in\fF(\BR^n_T)$ given by the family
$$\qquad\qquad\biggl\{u\in\fF(\BR^n_T):\|u\|_{L^{\iy}_T(L^1_{\om})}<\iy,
\lim_{s\to t^-}\|u(\cdot,s)-u(\cdot,t)\|_{L^1_{\om}}=0\text{
$\forall t\in(-T,0]$}\biggr\}$$ with the norm
$\|u\|_{L^{\iy}_T(L^1_{\om})}=\ds\sup_{t\in(-T,0]}\|u(\cdot,t)\|_{L^1_{\om}},$
which is separable with respect to the topology given by the norm.

\item A mapping $\BI:\fF(\BR^n_T)\to\fF(\BR^n_T)$ given by
$u\mapsto\BI u$ is called a {\it nonlocal parabolic operator} if (a)
$\BI u(x,t)$ is well-defined for any $u\in C^2_x(x,t)\cap
L^{\iy}_T(L^1_{\om})$ and (b) $\BI u$ is continuous on
$\Om_{\tau}\subset\BR^n_T$, whenever $u\in C^2_x(\Om_{\tau})\cap
L^{\iy}_T(L^1_{\om})$, where $C^2_x(x,t)$ is the class of all
$u\in\fF$ whose second derivatives $D^2 u$ in space variables exist
at $(x,t)$ and $C^2_x(\Om_{\tau})$ denotes the class of all
$u\in\fF$ such that $u\in C^2_x(x,t)$ for any $(x,t)\in\Om_{\tau}$
and $\ds\sup_{(x,t)\in\Om_{\tau}}|D^2 u(x,t)|<\iy$. Such a nonlocal
operator $\BI$ is said to be {\it uniformly elliptic} with respect
to a class $\fL$ of linear integro-differential operators if
\begin{equation}\bM^-_{\fL}v(x,t)\le
\BI(u+v)(x,t)-\BI u(x,t)\le\bM^+_{\fL}v(x,t)
\end{equation} where
$\bM^-_{\fL}v(x,t):=\inf_{\rL\in\fL}\rL v(x,t)$ and
$\bM^+_{\fL}v(x,t):=\sup_{\rL\in\fL}\rL v(x,t)$.

\item For $u\in C(Q_r)$, we define
$\|u\|_{C(Q_r)}=\sup_{(x,t)\in Q_r}|u(x,t)|$. For $\ap\in(0,1]$ and
$\sm\in(0,2)$, we define the {\it parabolic $\ap^{th}$ H\"older
seminorm} of $u$ by
$$\qquad\quad[u]_{C^{\ap}(Q_r)}=\sup_{(x,t),(y,s)\in
Q_r}\f{|u(x,t)-u(y,s)|}{(|x-y|^{\sm}+|t-s|)^{\ap/\sm}}.$$ In
particular, if $0<\ap/\sm<1$, then we define the norm
\begin{equation}\begin{split}\quad\|u\|_{C^{\sm+\ap}(Q_r)}&=\|u\|_{C(Q_r)}+\|\pa_t
u\|_{C(Q_r)}+\|(-\Delta)^{\sm/2}u\|_{C(Q_r)}\\
&+\|(Du)\mathbbm{1}_{[1,2)}(\sm)\|_{C(Q_r)}+[\pa_t
u]_{C^{\ap}(Q_r)}+[(-\Delta)^{\sm/2}u]_{C^{\ap}(Q_r)}.
\end{split}\end{equation}

\item For $a,b\in\BR$, we denote by $a\vee b=\max\{a,b\}$ and
$a\wedge b=\min\{a,b\}$.

\item For a multiindex $\bt=(\bt_1,\cdots,\bt_n)\in(\BN\cup\{0\})^n$,
we denote by $|\bt|=\sum_{i=1}^n\bt_i$.

\item Throughout this paper, {\it let $\e\in(0,1)$ be a fixed sufficiently
small positive number}.

\item For two quantities $a$ and $b$, we write $a\lesssim b$ (resp.
$a\gtrsim b$) if there is a universal constant $C>0$ ({\it depending
only on $\ld,\Ld,n,\e,\sm_0$ and the constants in $(1.3)$, $(1.4)$
and $(2.2)$, but not on $\sm$}) such that $a\le C\,b$ (resp. $b\le
C\,a$).

\item For $Q_r$, we denote by $C^2(Q_r)=C^2_x(Q_r)\cap C^1_t(Q_r)$ the class of functions $u\in\fF(\BR^n)$ which is $C^2$ in space and $C^1$ in time on $Q_r$.

\item For $(z,s)\in\BR^n_T$ and $u\in\fF(\BR^n_T)$, we denote the translation
operators $\btau_z$, $\btau^s$ and $\btau_z^s$ by $\btau_z
u(x,t)=u(x+z,t)$, $\btau^s u(x,t)=u(x,t+s)$ and $\btau_z^s
u(x,t)=u(x+z,t+s)$, respectively.

\item Let $f:\BR^n\times I\to\BR$ be a continuous function and let
$J:=(a,b]\subset I:=(-T,0]$. Then a function $u:\BR^n\times I\to\BR$
being upper $($lower$)$ semicontinuous on $\overline\Omega\times J$
is said to be a {\it viscosity subsolution} $($res. {\it viscosity
supersolution}$)$ of an equation $\BI u-\pa_t u=f$ on $\Omega\times
J$ and we write $\BI u-\pa_t u\ge f$ $($res. $\BI u-\pa_t u\le f$$)$
on $\Omega\times J$ in the viscosity sense, if for any
$(x,t)\in\Omega\times J$ there is a neighborhood
$Q_r(x,t)\subset\Omega\times J$ of $(x,t)$ such that $\BI
v(x,t)-\pa_t\vp(x,t)$ is well-defined and $\BI
v(x,t)-\pa_t\vp(x,t)\ge f(x,t)$ $($res. $\BI v(x,t)-\pa_t\vp(x,t)\le
f(x,t)$$)$ for
$v=\vp\mathbbm{1}_{Q_r(x,t)}+u\mathbbm{1}_{Q_r^c(x,t)}$ whenever
$\vp\in\rC^2(Q_r(x,t))$ with $\vp(x,t)=u(x,t)$ and $\vp>u$
$($$\vp<u$$)$ on $Q_r(x,t)\s\{(x,t)\}$ exists. Here, we denote such
a function $\vp$ by $\vp\in\rC^2_{\Om\times J}(u;x,t)^+$ $($res.
$\vp\in\rC^2_{\Om\times J}(u;x,t)^-$$)$. Also a function $u$ is
called as a {\it viscosity solution} if it is both a viscosity
subsolution and a viscosity supersolution to $\BI u-\pa_t u=f$ on
$\Omega\times J$ (see \cite{KL3, KL4}).

\item We say that $u\in\fF(\BR^n_T)$ is continuous at a point $(x,t)\in\pa_p Q_r$, if for any $\vep>0$ there exists some $\dt=\dt(\vep)>0$ such that
$|u(y,s)-u(x,t)|<\vep$ whenever $(y,s)\in\BR^n_T\s Q_r$ and $(|y-x|^{\sm}+|s-t|)^{1/\sm}<\dt$. If $u$ is continuous at every points in $\pa_p Q_r$, then we say that  $u$ is continuous in $\pa_p Q_r$ and we write $u\in C(\pa_p Q_r)$.

\end{itemize}

\,\,We shall now state the main theorem. The following
$\rC^{\sigma+\alpha}$-estimate for nonlocal parabolic concave
equation for $\sigma+\alpha\ge 2$ and $\sigma\in(1,2)$ makes it
possible to recover the well-known Evans-Krylov estimate as
$\sigma\rightarrow 2^-$. If $\sigma+\alpha<2$, then
$\rC^{\sigma+\alpha}$-estimate is covered by
$\rC^{1,\beta}$-estimate in \cite{KL3}. Our proof of the main
theorem is based on the nonlocal elliptic results of Silvestre and
Caffarelli \cite{CS1} and the regularity results on nonlocal
parabolic equations \cite{KL3, KL4}. Recently, we learned that Chang-Lara and Kriventsov \cite{CK} obtained some results for rough kernels under mild assumptions on the boundary data which is related with ours.

\begin{thm} Let $u\in L^{\iy}_T(L^1_\om)\cap C(\pa_p Q_2)$ be a viscosity solution
of the concave equation
\begin{equation*}\bI u-\pa_t u=0\,\text{ in
$Q_2$,}\end{equation*} where $\bI$ is defined on $\fL_2(\sm)$ for
$\sm\in(\sm_0,2)$ with $\sm_0\in(0,2)$ as in $(1.5)$. Then there
exists a constant $\ap\in(0,\f{1}{4}\wedge\sm_0\wedge|\sm_0-1|)$ such that
$$\|u\|_{C^{\sm+\ap}(Q_{1/2})}\lesssim\|u\|_{L^{\iy}_T(L^1_\om)}.$$
\end{thm}

\noindent{\bf Remark.} (i) As mentioned above, given any
$\sm_0\in(1,2)$ very close to $1$, it suffices to prove this theorem
only for $\sm+\ap\ge 2$ and $\sm\in[\sm_0,2)$.

(ii) In fact, from p.1569 of \cite{KL3} and (i) we could select such
$\ap>0$ so that $\ap\in(0,\f{1}{4}\wedge\sm_0\wedge|\sm_0-1|)$ in
the above. This implies that $0<\ap<2+\ap-\sm<1$.

\section{Parabolic interpolation inequalities}

Let $u\in C(Q_r)$. For $0<\ap\le 1$ and $\sm\in(0,2)$, we define the
{\it $\ap^{th}$ H\"older seminorms} of $u$ in the space and time
variable, respectively;

(i) $\ds
[u]_{C^{\ap}_x(Q_r)}=\sup_{t\in(-r^{\sm},0]}\,\sup_{(x,t),(y,t)\in
Q_r}\f{|u(x,t)-u(y,t)|}{|x-y|^{\ap}},$

(ii) $\ds [u]_{C^{\ap}_t(Q_r)}=\sup_{x\in B_r}\,\sup_{(x,t),(x,s)\in
Q_r}\f{|u(x,t)-u(x,s)|}{|t-s|^{\ap}}.$

\noindent If $\,0<\ap/\sm\le 1$, then it is easy to check that the
seminorms
$[\,\cdot\,]_{C^{\ap}_x(Q_r)}+[\,\cdot\,]_{C^{\f{\ap}{\sm}}_t(Q_r)}$
and $[\,\cdot\,]_{C^{\ap}(Q_r)}$ are equivalent.

We furnish an useful parabolic interpolation inequalities which
simplify the proof of our main.

\begin{thm} If $u\in L^{\iy}_T(L^1_{\om})\cap C(\pa_p Q_2)$ is a viscosity solution
of the concave equation
\begin{equation*}\bI u-\pa_t u=0\,\text{ in
$Q_2$,}\end{equation*} where $\bI$ is defined on $\fL_2(\sm)$ for
$\sm\in(\sm_0,2)$ with $\sm_0\in(1,2)$ very close to $1$, then there
exists a constant $\ap\in(0,1)$ with $\sm+\ap\ge
2$ such that
$$\|u\|_{C(Q_r)}\lesssim\|u\|_{L^{\iy}_T(L^1_\om)}\,\,\,\text{ and
}\,\,\, \|(-\Delta)^{\sm/2}u\|_{C(Q_r)}\vee\,\|\pa_t u\|_{C(Q_r)}\lesssim
\|u\|_{L^{\iy}_T(L^1_\om)}$$

\noindent for any $r\in(0,2)$.
\end{thm}

\pf By rescaling, the first and second inequalities can be shown as
in Theorem 5.1 and Corollary 7.4 below, respectively (also, refer to
\cite{KL4}). \qed

\noindent{\bf Remark.} (i) The main idea for the proof of the first
inequality comes from that of parabolic Harnack inequality, and so
it still holds without the concavity of the equation (see
\cite{KL3}).

(ii) Since $0<\ap/\sm<1$, this theorem and (1.8) imply that we have
only to control the seminorms $[\pa_t u]_{C^{\ap}(Q_r)}$ and
$[(-\Delta)^{\sm/2}u]_{C^{\ap}(Q_r)}$ in order to control the norm
$\|u\|_{C^{\sm+\ap}(Q_r)}$.

Next we give a fundamental lemma which facilitates the proof of
another type of parabolic interpolation inequalities.

\begin{lemma} If $\,u\in L^{\iy}_T(L^1_{\om})$ is a function with
$u(\cdot,t)\in C^k(B_r)$ for $t\in(-r^{\sm},0]$ and
$[D^{\bt}u]_{C_x^{\ap}(Q_r)}<\iy$ for some $\ap\in(0,1)$, then for
each $t\in(-r^{\sm},0]$ and multiindex $\bt$ with $|\bt|=k\in\BN$,
there exists some $z^t_0\in B_r$ $($depending on $t$$)$ such that
$$\bigl|D^{\bt}u(z^t_0,t)\bigr|\le
\bigl(\f{3r}{2}\bigr)^{\ap}\,[D^{\bt}u]_{C_x^{\ap}(Q_r)}+\f{2(4k)^k}{\om(B_{r/2})\,r^k}\|u\|_{L^{\iy}_T(L^1_{\om})}.$$
\end{lemma}

\pf Take $h=\f{r}{2k}$ and any multiindex $\bt$ with $|\bt|=k$. For
$(y,t)\in B_{r/2}\times(-T,0]$, we consider the finite difference
operator $\,D^{\bt}_h u(y,t)=D^{\bt_1}_{h,1}\,D^{\bt_2}_{h,2}\cdots
D^{\bt_n}_{h,n}u(y,t)$ where
$$D_{h,i}u(y,t)=\f{1}{h}\,[u(y+h e_i,t)-u(y,t)]$$
for a standard basis $\{e_1,\cdots,e_n\}$ of $\BR^n$. For
$i=1,\cdots,n$, we observe that
\begin{equation}D^{\bt_i}_{h,i}u(y,t)=\f{1}{h^{\bt_i}}\sum_{s=0}^{\bt_i}(-1)^{s}\f{\bt_i
!}{(\bt_i-s)!\,s !}\,u\bigl(y+(\bt_i-s)h e_i,t\bigr).
\end{equation}
By the mean value theorem, we see that there are some $z^t_1\in
B_h(y)$ and $z^t_2\in B_{2h}(y)$ such that
$$D_{h,i}D_{h,j}u(y,t)=\pa_{y_i}[D_{h,j}u](z^t_1,t)=D_{h,j}(\pa_{y_i}u)(z^t_1,t)
=\pa_{y_i y_j}u(z^t_2,t).$$ This implies that $D^{\bt}_h
u(y,t)=D^{\bt}u(z^t_y,t)$ for some $z^t_y\in B_{r/2}(y)$. Thus it
follows from this and (2.1) that
\begin{equation*}\begin{split}
\om(B_{r/2})\bigl|D^{\bt}u(z^t_0,t)\bigr|&\le\biggl|\om(B_{r/2})D^{\bt}u(z^t_0,t)-\int_{\BR^n}D_h^{\bt}u(y,t)\,\om(y)\,dy\biggr|
+\f{2^k}{h^k}\|u\|_{L^{\iy}_T(L^1_{\om})}\\
&\le\int_{B_{r/2}}\bigl|D^{\bt}u(z^t_0,t)-D^{\bt}u(z^t_y,t)\bigr|\,\om(y)\,dy+\f{2^{k+1}}{h^k}\|u\|_{L^{\iy}_T(L^1_{\om})}\\
&\le[D^{\bt}u]_{C_x^{\ap}(Q_r)}
\bigl(\f{3r}{2}\bigr)^{\ap}\,\om(B_{r/2})+\f{2(4k)^k}{r^k}\|u\|_{L^{\iy}_T(L^1_{\om})}.
\end{split}\end{equation*}
Therefore, this completes the proof. \qed

\begin{thm} If $\,u\in L^{\iy}_T(L^1_{\om})$ is a function such that
$u(\cdot,t)\in C^k(B_r)$ for each $t\in(-r^{\sm},0]$ and
$[D^{\bt}u]_{C_x^{\ap}(Q_r)}<\iy$ for some $\ap\in(0,1)$, then we
have that
$$\bigl\|D^{\bt}u\bigr\|_{C(Q_r)}\le
2\,\bigl(\f{3r}{2}\bigr)^{\ap}\,[D^{\bt}u]_{C_x^{\ap}(Q_r)}+\f{2(4k)^k}{\om(B_{r/2})\,r^k}\|u\|_{L^{\iy}_T(L^1_{\om})}$$
for any multiindex $\bt$ with $|\bt|=k\in\BN$.
\end{thm}

\pf From Lemma 2.2, for any $(x,t)\in Q_r$ we obtain that
\begin{equation*}\begin{split}\bigl|D^{\bt}u(x,t)\bigr|&\le\bigl|D^{\bt}u(x,t)-D^{\bt}u(z^t_0,t)\bigr|
+\bigl|D^{\bt}u(z^t_0,t)\bigr|\\
&\le 2\,[D^{\bt}u]_{C_x^{\ap}(Q_r)}
\bigl(\f{3r}{2}\bigr)^{\ap}+\f{2(4k)^k}{\om(B_{r/2})\,r^k}\|u\|_{L^{\iy}_T(L^1_{\om})}.
\end{split}\end{equation*}
Hence we can have the required inequality. \qed

In order to understand the {\bf parabolic H\"older spaces}
$C^{k,\gm}(Q_r)$ with $k\in\BN$ and $\gm\in(0,1)$, we define the
H\"older spaces $C^{k,\gm}_x(Q_r)$ and $C^{k,\gm}_t(Q_r)$ in the
space and time variable, respectively. For $u\in C(Q_r)$, we define
the norms
\begin{equation*}\begin{split}\|u\|_{C^{k,\gm}_x(Q_r)}&=\|u\|_{C(Q_r)}
+\sum_{i=1}^k\|D^i u\|_{C(Q_r)}+[D^k u]_{C^{\gm}_x(Q_r)},\\
\|u\|_{C^{k,\gm}_t(Q_r)}&=\|u\|_{C(Q_r)}+\sum_{i=1}^k\|\pa_t^i
u\|_{C(Q_r)} +[\pa_t^k
u]_{C^{\gm}_t(Q_r)},\end{split}\end{equation*} where $\|D^i
u\|_{C(Q_r)}=\sum_{|\bt|=i}\|D^{\bt}u\|_{C(Q_r)}$ and $[D^k
u]_{C^{\gm}_x(Q_r)}=\sum_{|\bt|=k}[D^{\bt}u]_{C^{\gm}_x(Q_r)}$ for
$i,k\in\BN$. And we denote by
$C_x^{k,\gm}(Q_r)=\{u\in\fF(\BR^n_T):\|u\|_{C^{k,\gm}_x(Q_r)}<\iy\}$
and
$C_t^{k,\gm}(Q_r)=\{u\in\fF(\BR^n_T):\|u\|_{C^{k,\gm}_t(Q_r)}<\iy\}$.

If $\sm\in[\sm_0,2)$ for $\sm_0\in(1,2)$ and $\ap\in(0,\sm_0-1)$,
then $0<\ap<2+\ap-\sm<1$ and
$$\f{2+\ap-\sm}{\sm}+1=\f{2+\ap}{\sm}.$$
Then we define the {\bf parabolic H\"older space} $C^{2,\ap}(Q_r)$
endowed with the norm
\begin{equation*}\begin{split}\|u\|_{C^{2,\ap}(Q_r)}&=\|u\|_{C(Q_r)}
+\sum_{i=1}^2\|D^i u\|_{C(Q_r)}+\|\pa_t u\|_{C(Q_r)}\\
&\qquad\qquad+[D^2 u]_{C^{\ap}(Q_r)}+[\pa_t u]_{C^{2+\ap-\sm}(Q_r)}.
\end{split}\end{equation*}

In the same case as the above, we can learn from Theorem 2.1 and
Theorem 2.3 that the estimates on the norm $\|u\|_{C^{2,\ap}(Q_r)}$
must be controlled by those on the seminorms $[\pa_t
u]_{C^{2+\ap-\sm}(Q_r)}\sim [\pa_t u]_{C_x^{2+\ap-\sm}(Q_r)}+[\pa_t
u]_{C_t^{\f{2+\ap-\sm}{\sm}}(Q_r)}$ and $\ds[D^2
u]_{C^{\ap}(Q_r)}\sim[D^2 u]_{C_x^{\ap}(Q_r)} +[D^2
u]_{C_t^{\f{\ap}{\sm}}(Q_r)}$. Similarly, the other parabolic
H\"older spaces can be defined along this line.

\begin{lemma} Let $\sm\in[\sm_0,2)$ for $\sm_0\in(1,2)$ and $\ap\in(0,\sm_0-1)$.
If $u\in L^{\iy}_T(L^1_{\om})$ is a function with $u(x,\cdot)\in
C^1(-r^{\sm},0]$ for $x\in B_r$ and $[\pa_t
u]_{C_t^{\f{2+\ap-\sm}{\sm}}(Q_r)}<\iy$, then we have that
\begin{equation*}\|\pa_t u\|_{C(Q_r)}\le r^{2+\ap-\sm}[\pa_t
u]_{C_t^{\f{2+\ap-\sm}{\sm}}(Q_r)}+\f{4}{r^{\sm}}\|u\|_{C(Q_r)}.
\end{equation*}
\end{lemma}

\pf Take any $r\in(0,2)$ and $(x,t)\in Q_r$. Then there is some
$t_0\in(-r^{\sm},0]$ such that $|t-t_0|=\f{1}{2}r^{\sm}$, and by the
mean value theorem, there is some $t_0^x$ between $t$ and $t_0$ such
that $u(x,t_0)-u(x,t)=\f{1}{2}r^{\sm}\,\pa_t u(x,t_0^x)$. Thus we
have the estimate
\begin{equation*}\begin{split}
\f{1}{2}\,r^{\sm}\,|\pa_t u(x,t)|&\le\biggl|\f{1}{2}\,r^{\sm}\,\pa_t
u(x,t)-\bigl(u(x,t_0)-u(x,t)\bigr)\biggr|+2\,\|u\|_{C(Q_r)}\\
&=\f{1}{2}\,r^{\sm}\bigl|\pa_t u(x,t)-\pa_t
u(x,t_0^x)\bigr|+2\,\|u\|_{C(Q_r)}\\
&\le\f{1}{2}\,r^{2+\ap}[\pa_t
u]_{C_t^{\f{2+\ap-\sm}{\sm}}(Q_r)}+2\,\|u\|_{C(Q_r)}.
\end{split}\end{equation*}
Hence this implies the required inequality. \qed

\begin{lemma} Let $\sm\in[\sm_0,2)$ for $\sm_0\in(1,2)$, and
let $u\in L^{\iy}_T(L^1_{\om})$ be a viscosity solution of
the equation
$$\bI u-\pa_t u=0\,\,\text{ in $Q_2$}$$ where $\bI$ is defined on $\fL_2(\sm)$. If $u\in $, then we have the
estimates
\begin{equation*}\begin{split}[D^2 u]_{C_t^{\f{\ap}{\sm}}(Q_r)}&\lesssim
\|D^2 u\|_{C(Q_r)}+\|u\|_{L^{\iy}_T(L^1_{\om})},\\
[\pa_t u]_{C_x^{2+\ap-\sm}(Q_r)}&\lesssim\|u\|_{L^{\iy}_T(L^1_{\om})}\end{split}\end{equation*}
for any $r\in(0,1)$.
\end{lemma}

\pf Take any $r\in(0,2)$ and $(x,t)\in Q_r$. We note that
$0<\ap<2+\ap-\sm<1$. For $h$ with $|h|<\ep$, we consider the
difference quotients in the $x$-direction
$$u^h(x,t)=\f{u(x+h,t)-u(x,t)}{|h|}.$$ Write $u^h=u^h_1+u^h_2$ where
$u^h_1=u^h\mathbbm{1}_{Q_r}$. By Theorem 2.4 \cite{KL3}, we have
that $\bM^+_{\fL_2}u^h-\pa_t u^h\ge 0$ and $\bM^-_{\fL_2}u^h-\pa_t
u^h\le 0$ on $Q_r$. Since $\pa_t u^h_2\equiv 0$ in $Q_r$, it follows
from the uniform ellipticity (1.7) of $\bM^+_{\fL_2}$ and
$\bM^-_{\fL_2}$ with respect to $\fL_2$ that
$$\bM^+_{\fL_0}u^h_1-\pa_t u^h_1\ge -\bM^+_{\fL_2}u^h_2\,\,\text{ and
}\,\,\bM^-_{\fL_0}u^h_1-\pa_t u^h_1\le -\bM^-_{\fL_2}u^h_2\,\,\text{
in $Q_r$.}$$ Then it is easy to show that
$|\bM^+_{\fL_2}u^h_2|\vee|\bM^-_{\fL_2}u^h_2|\lesssim
\|u\|_{L^{\iy}_T(L^1_{\om})}$ in $Q_r$ for a universal constant
$c>0$. So we have that $$\bM^+_{\fL_0}u^h_1-\pa_t u^h_1\gtrsim
-\|u\|_{L^{\iy}_T(L^1_{\om})}\,\,\text{ and
}\,\,\bM^-_{\fL_0}u^h_1-\pa_t u^h_1\lesssim
\|u\|_{L^{\iy}_T(L^1_{\om})}\,\,\text{ in $Q_r$.}$$ We now
consider another difference quotients in the $x$-direction
$$w^h(x,t)=\f{u^h_1(x+h,t)-u^h_1(x,t)}{|h|}.$$
Applying Theorem 2.4 \cite{KL3} again, we obtain that
$$\bM^+_{\fL_0}w^h-\pa_t w^h\gtrsim
-\|u\|_{L^{\iy}_T(L^1_{\om})}\,\text{ and
}\,\,\bM^-_{\fL_0}w^h-\pa_t w^h\lesssim
\|u\|_{L^{\iy}_T(L^1_{\om})}\,\text{ in $Q_r$.}$$ From the
H\"older estimate(Theorem 3.4) in \cite{KL4}, we get the estimate
$$[w^h]_{C^{\f{\ap}{\sm}}_t(Q_r)}\le[w^h]_{C^{\ap}(Q_r)}\lesssim
\|w^h\|_{C(Q_r)}+\|w^h\|_{L^{\iy}_T(L^1_{\om})}+\|u\|_{L^{\iy}_T(L^1_{\om})}.$$
By the mean value theorem, we easily have that $\|w^h\|_{C(Q_{r_{\ep}})}\le\|D^2
u\|_{C(Q_{r_{\ep}})}$. Since $|D\om(y,s)|+|D^2\om(y,s)|\lesssim
\om(y)$, it follows from the integration by parts that
$\,\|w^h\|_{L^{\iy}_T(L^1_{\om})}\le\|u\|_{L^{\iy}_T(L^1_{\om})}.$
Thus we obtain that
$$[w^h]_{C^{\f{\ap}{\sm}}_t(Q_r)}\lesssim
\|D^2 u\|_{C(Q_r)}+\|u\|_{L^{\iy}_T(L^1_{\om})}.$$
Taking the limit $|h|\to 0$, the first inequality can be obtained.

Take any $(x,t)\in Q_r$. Then it follows from the uniform ellipticity that
\begin{equation}\begin{split}\bM^-_2(\btau^t_x u-\btau^t u)(0,0)&\le\pa_t u(x,t)-\pa_t u(0,t)\\&=\bI u(x,t)-\bI u(0,t)
\le\bM^+_2(\btau^t_x u-\btau^t u)(0,0)
\end{split}\end{equation}
Let $\vp\in C^{\iy}_c(\BR^n)$ be a function satisfying that $\vp=1$ in
$B_1$, $\vp=0$ in $\BR^n\s B_{3/2}$ and $0\le\vp\le 1$ in $\BR^n$, and take any $\rL\in\fL_2$. Then by the change of variable, the mean value theorem and (1.3) we have that
\begin{equation}\begin{split}\rL(\btau_x^t u-\btau^t u)(0,0)&=\int_{\BR^n}\bigl[\mu_t(u,x,y)-\mu_t(u,0,y)\bigr]\vp(y)K(y)\,dy\\
&+\int_{\BR^n}\bigl[\mu_t(u,x,y)-\mu_t(u,0,y)\bigr](1-\vp(y))K(y)\,dy\\
&\lesssim\vp^+ u(x,0)+\|u\|_{L^{\iy}_T(L^1_{\om})}\,|x|,
\end{split}\end{equation}
where $$\vp^+ u(x,0)=\sup_{t\in(-T,0]}\sup_{K\in\cK_2}\int_{\BR^n}\bigl[\mu_t(u,x,y)-\mu_t(u,0,y)\bigr]\vp(y)K(y)\,dy.
$$
Similarly we can obtain that
\begin{equation}\begin{split}\rL(\btau_x^t u-\btau^t u)(0,0)\gtrsim\vp^- u(x,0)-\|u\|_{L^{\iy}_T(L^1_{\om})}\,|x|,
\end{split}\end{equation}
where $$\vp^- u(x,0)=\inf_{t\in(-T,0]}\inf_{K\in\cK_2}\int_{\BR^n}\bigl[\mu_t(u,x,y)-\mu_t(u,0,y)\bigr]\vp(y)K(y)\,dy.
$$ The estimates (2.2), (2.3) and (2.4) imply that
\begin{equation}\begin{split}\vp^- u(x,0)-\|u\|_{L^{\iy}_T(L^1_{\om})}\,|x|&\lesssim\bM^-_2(\btau_x^t u-\btau^t u)(0,0)\\&\le\pa_t u(x,t)-\pa_t u(0,t)\\
&\le\bM^+_2(\btau_x^t u-\btau^t u)(0,0)\\
&\lesssim\vp^+ u(x,0)+\|u\|_{L^{\iy}_T(L^1_{\om})}\,|x|.
\end{split}\end{equation}
Applying the method in Lemma 9.2 \cite{CS1} with (2.5), we have that
$$|\vp^- u(x,0)|\vee|\vp^+ u(x,0)|\lesssim\|u\|_{L^{\iy}_T(L^1_{\om})}\,|x|^{\bt}$$
for some $\bt\in(0,1)$. Here, without loss of generality, we may assume that $\bt=2+\ap-\sm$ by applying a standard telescopic argument \cite{CC}.
Hence the second inequality can be achieved from a standard translation argument. Therefore we complete the proof. \qed

We now consider the class $\fL_*$ of operators $L$ with kernels
$K\in\cK_*$ satisfying (1.2) such that there are some $\vr_0>0$ and a constant $C>0$ such that
\begin{equation}|\n K(y)|\le
C\,\om(y)\,\,\text{ for any $y\in\BR^n\s B_{\vr_0}.$ }\end{equation}
We note that $\fL_1$ is the largest scale invariant class contained
in the class $\fL_*$.

\begin{thm} Let $\sm\in[\sm_0,2)$ for some $\sm_0\in(1,2)$.
Then there is some $\vr_0>0$ $($depending on $\ld,\Ld,\sm_0$ and
$n$$)$ so that if $\bI$ is a nonlocal, translation-invariant and
uniformly elliptic operator with respect to $\fL_*$ and $u\in
L^{\iy}_T(L^1_{\om})$ satisfies the equation
$$\bI u-\pa_t u=0\,\,\text{ in $Q_2$,}$$ then
there is some $\ap>0$ such that
$$\|D u\|_{C_t^{\f{\ap}{\sm}}(Q_r)}\lesssim
\|D u\|_{C(Q_r)}+\|u\|_{L^{\iy}_T(L^1_{\om})}$$ for
any $r\in(0,2)$.\end{thm}

\pf We proceed the proof by applying Theorem 3.4 \cite{KL4} to the
difference quotients in the $x$-direction
$$w^h(x,t)=\f{u(x+h,t)-u(x,t)}{|h|}.$$
Take any $r\in(0,2)$. Then we write $w^h=w_1^h+w_2^h$ where
$w_1^h=w^h\mathbbm{1}_{Q_r}$. From Theorem 2.4 \cite{KL3}, we have
that $\bM^+_{\fL^*}w^h-\pa_t w^h\ge 0$ and $\bM^-_{\fL^*}w^h-\pa_t
w^h\le 0$ in $Q_r$. Because $\pa_t w_2^h\equiv 0$ in $Q_r$, it
follows from the uniform ellipticity with respect to $\fL^*$ that we
get that
\begin{equation*}\begin{split} \bM^+_{\fL_0}w_1^h-\pa_t w_1^h&\ge
\bM^+_{\fL_*}w_1^h-\pa_t w_1^h\ge\bM^+_{\fL_*} w^h-\bM^+_{\fL_*}
w_2^h-\pa_t w^h\ge-\bM^+_{\fL_*}w_2^h\,\,\text{ in $Q_r$},\\
\bM^-_{\fL_0}w_1^h-\pa_t w_1^h&\le\bM^-_{\fL_*}w_1^h-\pa_t
w_1^h\le\bM^-_{\fL_*} w^h-\bM^-_{\fL_*} w_2^h-\pa_t
w^h\le-\bM^-_{\fL_*}w_2^h\,\,\text{ in
$Q_r$}.\end{split}\end{equation*} If we can show that
$|\bM^+_{\fL_*}w_2^h|\vee|\bM^-_{\fL_*}w_2^h|\lesssim
\|u\|_{L^{\iy}_T(L^1_{\om})}$ in $Q_r$, then we have that
$$\bM^+_{\fL_0}w_1^h-\pa_t w_1^h\gtrsim-\|u\|_{L^{\iy}_T(L^1_{\om})}\text{ and
}\bM^-_{\fL_0}w_1^h-\pa_t w_1^h\lesssim
\|u\|_{L^{\iy}_T(L^1_{\om})}\text{ in $Q_r$}$$ for $h$ with a
sufficiently small $|h|$. Indeed, by using (2.6), it can be obtained
from the fact that
\begin{equation*}\begin{split}&\int_{\BR^n\s
B_{\rho}}|u(x+y,t)|\f{|K(x,y,t)-K(x,y-h,t)|}{|h|}\,dy\\
&\qquad+\int_{\BR^n\s B_{\rho}}|u(x+y+h,t)|K(x,y,t)\,dy\lesssim
\|u\|_{L^{\iy}_T(L^1_{\om})}\end{split}\end{equation*} for some
$\rho>0$. Hence $w_1^h$ admits the H\"older estimate(Theorem 3.4 \cite{KL4})
on $Q_r$, and thus applying the mean value theorem and
integration by parts with (2.6) gives the estimate
\begin{equation*}\begin{split}\|w_1^h\|_{C_t^{\f{\ap}{\sm}}(Q_r)}
&\le\|D u\|_{C(Q_r)}+\|u\|_{L^{\iy}_T(L^1_{\om})}.
\end{split}\end{equation*}
Finally, taking the limit $|h|\to 0$, we obtain the required result.
\qed

{\bf Remark.} In order to show Theorem 1.1, we learned from the
interpolation results obtained in this section that the norm
$\|u\|_{C^{2,\ap}(Q_r)}$ of viscosity solutions $u\in
L^{\iy}_T(L^1_{\om})$ of the equation
$$\bI u-\pa_t u=0\,\,\text{ in $Q_r$}$$
is controlled by only two
seminorms $[\pa_t u]_{C_t^{\f{2+\ap-\sm}{\sm}}(Q_r)}$ and $[D^2
u]_{C_x^{\ap}(Q_r)}$, and so only two norms
$\|u\|_{C^{2,\ap}_x(Q_r)}$ and
$\|u\|_{C^{1,\f{2+\ap-\sm}{\sm}}_t(Q_r)}$.

\section{Approximation of solutions and average of subsolutions}

In the first part of this section, we show that any viscosity
solution of (1.5) can be approximated by $C^{2,\ap}$-functions
solving an approximate equation with the same shape as (1.5), by
using a standard regularization argument. This useful result makes
it possible to extend an estimate on $C^{2,\ap}$-solutions to the
estimate on viscosity solutions by passing to the limit process.

Let $\Om$ be a bounded domain in $\BR^n$ and
$\Om_{\tau}=\Om\times(-\tau,0]$ for $\tau\in(0,T)$. Then we say that
a function $u:\BR^n_T\to\BR$ is in $C^{1,1}_x(\Om_{\tau})$, if there
is a constant $C_0>0$ (independent of $(x,t)$ and $(y,t)$) such that
\begin{equation}|u(y,t)-u(x,t)-(y-x)\cdot\n u(x,t)|\le C_0|y-x|^2\end{equation} for all
$(x,t),(y,t)\in\Om_{\tau}$. Here we denote by the norm
$\|u\|_{C^{1,1}_x(\Om_{\tau})}$ the smallest $C_0$ satisfying (3.1).

The following definitions are the parabolic version corresponding to
the elliptic case in \cite{CS1} (see also \cite{KL4}).

\begin{definition} For a nonlocal parabolic operator $\rI$ and $\tau\in(0,T)$, we define
$\|\rI\|$ in $\Om_{\tau}$ with respect to a weight $\om$ as
$$\|\rI\|=\sup_{(y,s)\in\Om_{\tau}}\sup_{u\in\cF^M_{y,s}}\f{|\rI u(y,s)|}{1+\|u\|_{L^{\iy}_T(L^1_{\om})}+\|u\|_{C_x^{1,1}(Q_1(y,s))}}$$
where $\cF^M_{y,s}=\{u\in\fF(\BR^n_T)\cap
 C^2_x(y,s):\|u\|_{L^{\iy}_T(L^1_{\om})}\vee\|u\|_{C_x^{1,1}(Q_1(y,s))}\le M\}$ for some $M>0$.
\end{definition}

For $K_{\bt}\in\fL_0$ and $\vep>0$, we consider the following
regularized kernels
$$K^{\vep}_{\bt}(y)=\vp_{\vep}(y)\f{\ld(2-\sm)}{|y|^{n+\sm}}+(1-\vp_{\vep}(y))K_{\bt}(y)$$
where $\vp\in C^{\iy}_c(\BR^n)$ is a function such that $0\le\vp\le
1$ in $\BR^n$, $\vp=0$ in $\BR^n\s B_2$ and $\vp=1$ in $B_1$, and
$\vp_{\vep}(y)=\vp(y/\vep)$. Then we define the corresponding
operator $\bI^{\vep}$ by
$$\bI^{\vep}v(x,t):=\inf_{\bt}\rL^{\vep}_{\bt}v(x,t):=\inf_{\bt}\int_{\BR^n}\mu_t(v,x,y)K_{\bt}^{\vep}(y)\,dy.$$
Under the parabolic topology, it is natural to consider the partial
derivative $\pa^-_t$ with respect to the past time defined by
$$\pa^-_t u(x,t)=\lim_{h\to 0^-}\f{u(x,t+h)-u(x,t)}{h}$$ for
$u\in\frak F$, if it exists.

\begin{lemma} Let $u\in L^{\iy}_T(L^1_{\om})\cap C(\pa_p Q_{1+\e})$ be a
viscosity solution of the nonlocal parabolic concave
equation
$$\bI u-\pa_t u=0\,\,\text{ in $Q_{1+\e}$, }$$
where every $\rL_{\bt}$ belong to the class $\fL_2(\sm)$ for $\sm\in(1,2)$.
Then there are some $\ap\in(0,1)$
and a sequence $\{u^{\vep}\}\subset C^{2,\ap}(Q_1)$ such that
$$\lim_{\vep\to 0}\,\sup_{Q_{1+\e}}|u^{\vep}-u|=0,\,\,\,\,
\lim_{\vep\to 0}\pa_t u^{\vep}=\pa_t u\text{ on
$B_1\times(-1,0),$}$$ $\lim_{\vep\to 0}\pa_t u^{\vep}(x,0)=\pa^-_t
u(x,0)$ for any $x\in B_1$ and
\begin{equation}\begin{cases}
\bI^{\vep}u^{\vep}-\pa_t u^{\vep}=0 &\text{ in $Q_{1+\e}$, }\\
u^{\vep}=u &\text{ in $\BR^n_T\s Q_{1+\e}.$}
\end{cases}\end{equation}
Moreover, we have that $\lim_{\vep\to 0}\|\bI^{\vep}-\bI\|=0$.
\end{lemma}

{\bf Remark.} Note that the condition $\lim_{\vep\to
0}\|\bI^{\vep}-\bI\|=0$ implies that $\rI^{\vep}$ converges weakly
to $\rI$ in $Q_{1+\e}$ as in \cite{KL4}.

\,\,\pf We observe that if $\rL_{\bt}\in\fL_2(\sm)$, then
$\rL_{\bt}^{\vep}\in\fL_2(\sm)$. For any $\vep\in(0,1)$, let $u^{\vep}$
be the viscosity solution of (3.2). Then it follows from Corollary 7.9
\cite{KL4} that $u^{\vep}\in C^{2,\ap}(Q_1)$ for some $\ap\in(0,1)$.

If $v\in \cF^M_{y,s}$ for $M>0$ and $(y,s)\in Q_1$, then
$\|v\|_{L^{\iy}_T(L^1_{\om})}\vee\|v\|_{C^{1,1}(Q_1(y,s))}\le M$ and
$v\in\fF\cap C^2_x(y,s)$, and so we have that
$$|v(x,t)-v(y,t)-(x-y)\cdot\n_x
v(y,t)|\le\|v\|_{C^{1,1}(Q_1(y,s))}|x-y|^2$$ for all $(x,t)\in
Q_1(y,s)$. Thus by simple computation, we obtain that
$$|\bI^{\vep}v(x,t)-\bI v(x,t)|\lesssim
\vep^{2-\sm},$$ so that $\|\bI^{\vep}-\bI\|\lesssim\vep^{2-\sm}\to
0$ as $\vep\to 0$ because $\sm\in(0,2)$. Thus by Lemma 5.8
\cite{KL4} we conclude that $u^{\vep}$ converges to $u$ uniformly in
$Q_{1+\e}$ as $\vep\to 0$.

For $\vep\in(0,1)$, $h\in(-1,1)$ and $(x,t)\in Q_1$, we set
$$g_{\vep,h}(x,t)=\f{u^{\vep}(x,t+h)-u^{\vep}(x,t)}{h}\,\,\text{ and
}\,\,g_h(x,t)=\f{u(x,t+h)-u(x,t)}{h}.$$ For every fixed $h\in(0,1)$,
it is easy to check that $g_{\vep,h}$ converges uniformly to $g_h$
on $Q_1$ as $\,\vep\to 0$, and moreover $g_{\vep,h}$ has a pointwise
limit $\pa_t u^{\vep}$ on $Q_1$ as $\,h\to 0$. Thus, by commutative
property of double limits, $g_h$ has a pointwise limit on $Q_1$ as
$h\to 0$, and moreover $$\pa_t u(x,t)=\lim_{h\to
0}g_h(x,t)=\lim_{\vep\to 0}\lim_{h\to
0}g_{\vep,h}(x,t)=\lim_{\vep\to 0}\pa_t u^{\vep}(x,t)$$ for any
$(x,t)\in B_1\times(-1,0)$ and $\lim_{\vep\to 0}\pa_t
u^{\vep}(x,0)=\pa^-_t u(x,0)$ for any $x\in B_1$. Hence we are done.
\qed

From Lemma 3.2 and Theorem 2.2 \cite{KL3}, we can easily derive the
following corollary which shall be useful in the final step of the
proof of the main theorem.

\begin{cor} If $u\in L^{\iy}_T(L^1_{\om})\cap C(\pa_p Q_{1+\e})$ be a
viscosity solution of the nonlocal parabolic concave equation
$$\bI u-\pa_t u=0\,\,\text{ in $Q_{1+\e}$, }$$
where every $\rL_{\bt}$ belong to $\fL_2(\sm)$ for $\sm\in(1,2)$,
then $\bI u-\pa_t u$ is well-defined on
$Q_{1+\e}$ in the classical sense and
\begin{equation*}\bI u(x,t)-\pa_t u(x,t)=0\,\,\text{ for any $(x,t)\in Q_{1+\e}$.}
\end{equation*}
\end{cor}

In the second part, we shall show that any average of viscosity
subsolutions to the nonlocal parabolic concave equation is a
viscosity subsolution to the same equation. This implies that the
convolution of the viscosity subsolution with a mollifier with
compact support is also a viscosity subsolution, which shall be very
useful in obtaining local uniform boundedness of linear operators in
Section 6.

\begin{lemma} If $u,v\in L^{\iy}_T(L^1_{\om})\cap C(\pa_p Q_1)$ be
viscosity subsolutions of the concave equations $\bI u-\pa_t u=0$
and $\bI v-\pa_t v=0$ in $Q_1$, then we have that
$$\bI\biggl(\f{u+v}{2}\biggr)-\pa_t\biggl(\f{u+v}{2}\biggr)\ge 0\,\,\text{ in $Q_1$ }$$
in the viscosity sense. In particular, if $\,u\in
L^{\iy}_T(L^1_{\om})\cap C(\pa_p Q_1)$ is a viscosity solution of the concave
equation $\bI u-\pa_t u=0$ in $Q_1$  and $\vp\in C^{\iy}_c(\BR^n)$
is a mollifier supported in a small ball $B_{\dt}$ such that $\vp\ge
0$ and $\|\vp\|_{L^1(\BR^n)}=1$, then $\bI(\vp*u)-\pa_t(\vp*u)\ge 0$
in $Q_1$ in the viscosity sense.
\end{lemma}

\rk Note that the convolution $\vp*u$ of $\vp$ and $u$ means
$$\vp*u(x,t)=\int_{\BR^n}\vp(x-y)u(y,t)\,dy=\int_{\BR^n}u(x-y,t)\vp(y)\,dy,\,\,x\in\BR^n,t\in(-T,0].$$

\pf We consider approximate equations $\bI^{\vep}u^{\vep}-\pa_t
u^{\vep}=0$ and $\bI^{\vep}v^{\vep}-\pa_t v^{\vep}=0$ in $Q_1$ with
boundary values as in (3.2). By Lemma 3.2, we see that
$u^{\vep},v^{\vep}\in C^2(Q_1)$ and $u^{\vep},v^{\vep}$ converges
uniformly to $u,v$ in $Q_1$, respectively. Thus the operators
$\rL^{\vep}_{\bt}u^{\vep},\rL^{\vep}_{\bt}v^{\vep}$ are well-defined
and continuous on $Q_1$. Now it follows from simple computation that
\begin{equation*}\begin{split}\bI^{\vep}\biggl(\f{u^{\vep}+v^{\vep}}{2}\biggr)-\pa_t\biggl(\f{u^{\vep}+v^{\vep}}{2}\biggr)
&\ge\f{\inf_{\bt}\rL_{\bt}u^{\vep}+\inf_{\bt}\rL_{\bt}v^{\vep}}{2}-\pa_t\biggl(\f{u^{\vep}+v^{\vep}}{2}\biggr)\\
&=\f{(\bI^{\vep}u^{\vep}-\pa_t u^{\vep})+(\bI^{\vep}v^{\vep}-\pa_t
v^{\vep})}{2}\ge 0\,\,\text{ in $\Om\times I$ }
\end{split}\end{equation*} in the viscosity sense. Since it is
obvious that $\lim_{\vep\to
0}\|u^{\vep}-u\|_{L^{\iy}_T(L^1_{\om})}=0$ and $\lim_{\vep\to
0}\|v^{\vep}-v\|_{L^{\iy}_T(L^1_{\om})}=0$, by Lemma 5.4 \cite{KL4}
and Lemma 3.2 we obtain the first required result. Finally, the
second part is a natural by-product of the first part we obtained
just before in the above. \qed

\section{Linear parabolic integro-differential equations}

In this section, we shall obtain regularity results for linear
parabolic integro-differential equations much better than those for
the nonlinear equations.

\begin{thm} Let $\rL$ be a linear integro-differential operator in
the class $\fL_1(\sm)$ for $\sm\in(\sm_0,2)$ with $\sm_0\in(1,2)$.
If $u\in L^{\iy}_T(L^1_{\om})\cap C(\pa_p Q_{1+\e})$ is a viscosity solution of $$\rL
u-\pa_t u=0\,\,\text{ in $Q_{1+\e},$ }$$ then $u\in C^{2,\ap}(Q_1)$,
and moreover there is some $\ap\in(0,1)$ such that
$$\|u\|_{C^{2,\ap}(Q_1)}\lesssim\|u\|_{L^{\iy}_T(L^1_{\om})}.$$
\end{thm}

\pf Applying Theorem 3.6 in \cite{KL4}, we see that there is a
constant $\ap\in(0,1)$ such that $u\in C^{1,\ap}(Q_1)$ and
\begin{equation}\|u\|_{C^{1,\ap}(Q_1)}\lesssim
\|u\|_{L^{\iy}_T(L^1_{\om})}.\end{equation}
We note that $\rL u_e(x,t)-\pa_t u_e(x,t)=0$ for $(x,t)\in Q_1$
where $u_e$ means the weak derivative of $u$ in the direction $e\in
S^{n-1}$. Also by (4.1), we note that $u_e$ coincides with the
strong type directional derivative of $u$ in the direction $e$ on
$Q_1$.

Next we show that $u_e\in L^{\iy}_T(L^1_{\om})$. For $(x,t)\in Q_1$,
we consider a function $w\in C^1_0(\BR^n)$ such that $w(y)=1$ for
$|y|<1/2$, $|w_e(y)|\le 1$ and $w(y)\ge 1$ for $1/2\le|y|<1$, and
$w(y)=K(y)$ for $|y|\ge 1$. Take any $(x,t)\in Q_1$. Then by
integration by parts, (1.3) and Theorem 2.1, we have that
\begin{equation}\begin{split}&\biggl|\int_{\BR^n}u_e(y,t)w(y)\,dy\biggr|
=\biggl|\int_{\BR^n}u(y,t)w_e(y)\,dy\biggr|\\
&\qquad\quad\le\int_{1/2\le|y|<1}|u(y,t)|\,dy+\int_{|y|\ge
1}|u(y,t)|\,|\la e,\n K(y)\ra|\,dy\\
&\qquad\quad\lesssim
\|u\|_{C(Q_{1+\e})}+\|u\|_{L^{\iy}_T(L^1_{\om})}\lesssim\|u\|_{L^{\iy}_T(L^1_{\om})}.
\end{split}\end{equation}
This implies that $(u_e)^+,(u_e)^-\in L^1(w\,dy)$. Thus we see that
$|u_e|\in L^1(w\,dy)$. Moreover we conclude that $u_e\in
L^{\iy}_T(L^1_{\om})$.

Then it follows from Theorem 3.6 \cite{KL4} that $u_e\in
C^{1,\ap}_x(Q_1)$. Thus we obtain that $u\in C^{2,\ap}_x(Q_1)$. Here
we note that we could choose some $\ap>0$ so that $\ap<\sm_0-1$ in
Theorem 3.4 \cite{KL4} (or Theorem 5.2 \cite{KL3}). Since $(2+\ap)/\sm>1$ for such $\ap>0$, we
see that
$$\f{2+\ap-\sm}{\sm}+1=\f{2+\ap}{\sm}$$ and $0<\ap<2+\ap-\sm<1$.
Since $0<2+\ap-\sm<1<1+\ap$, by (4.1) we can obtain that $u$ is
$C_t^{\f{2+\ap-\sm}{\sm}}$-H\"older continuous in $Q_1$. By applying the idea
of the proof of Theorem 7.8 \cite{KL4}, the
$C^{1,\f{2+\ap-\sm}{\sm}}_t$-regularity of $u$ can be achieved on
$Q_1$. Therefore by the final remark in Section 2, we conclude that
$u\in C^{2,\ap}(Q_1)$.\qed

Let $\fF(\BR^n_T)$ denote the family of all real-valued measurable
functions defined on $\BR^n_T$. Then we introduce a function space
$L^{\iy}_T(L^2_x)$ consisting of all $f\in\fF(\BR^n_T)$ satisfying
$$\sup_{t\in(-T,0]}\biggl(\int_{\BR^n}|f(x,t)|^2\,dx\biggr)^{\f{1}{2}}<\iy.$$

\begin{thm} If $\|\rL_0 u\|_{L^{\iy}_T(L^2_x)}<\iy$ for some
$\rL_0\in\fL_0(\sm)$ with $\sm\in(0,2)$ and $u\in\fF$, then we have
that
$$\sup_{\rL\in\fL_0(\sm)}\|\rL u\|_{L^{\iy}_T(L^2_x)}
\lesssim \inf_{\rL\in\fL_0(\sm)}\|\rL u\|_{L^{\iy}_T(L^2_x)}.$$
\end{thm}

\pf If we denote the Fourier transform of $u\in\fF(\BR^n_T)$ in
terms of space variable by $\widehat u(\xi,t)=\int_{\BR^n}e^{-i
x\cdot\xi}u(x,t)\,dx,$ then it follows from Plancherel's Theorem
that $$\widehat{\rL
u}(\xi,t)=\biggl(-\int_{\BR^n}2\bigl(1-\cos(y\cdot\xi)\bigr)K(y)\,dy\biggr)\widehat
u(\xi,t):=-m(\xi)\widehat u(\xi,t)$$ for any $\rL\in\fL_0(\sm)$. By
simple computation as in \cite{CS1}, we have that
$$\f{1}{c_0}|\xi|^{\sm}\le m(\xi)\le c_0|\xi|^{\sm}$$ for a
universal constant $c_0>0$ possibly depending on $\ld,\Ld$ and the
dimension $n$, but not depending on $t$. Applying standard harmonic
analysis, there is a universal constant $C>0$ possibly depending on
$\ld,\Ld$ and the dimension $n$, but not depending on $t$ such that
$$\sup_{t\in(-T,0]}\,\sup_{\|v(\cdot,t)\|_{L^2(\BR^n)}\neq 0}
\f{\|\rL_1\circ\rL_2^{-1}v(\cdot,t)\|_{L^2(\BR^n)}}{\|v(\cdot,t)\|_{L^2(\BR^n)}}
=\|m_1\,m_2^{-1}\|_{L^{\iy}(\BR^n)}<C<\iy$$ for any
$\rL_1,\rL_2\in\fL_0(\sm)$, where $m_1$ and $m_2^{-1}$ denote the
symbols of $\rL_1$ and  the inverse $\rL_2^{-1}$ of the operator
$\rL_2$, respectively. Hence this implies the required result. \qed

Let $s$ be a real number. Then the homogeneous mixed Sobolev space
$L^{\iy}_T({\dot{H}}^s_x)$ is defined as the function space of all
$f\in\fF(\BR^n_T)$ satisfying
$$\|f\|_{L^{\iy}_T({\dot{H}}^s_x)}
:=\sup_{t\in(-T,0]}\biggl(\int_{\BR^n}|\xi|^{2s}\bigl|\widehat
f(\xi,t)\bigr|^2\biggr)^{\f{1}{2}}<\iy.$$ For $p_{\sm}:=2n/(n-2\sm)$
with $\sm\in(0,2)$, we define a function space
$L^{\iy}_T(L^{p_{\sm}}_x)$ consisting of all $f\in\fF(\BR^n_T)$
satisfying
$$\|f\|_{L^{\iy}_T(L^{p_{\sm}}_x)}
:=\sup_{t\in(-T,0]}\biggl(\int_{\BR^n}\bigl|f(x,t)\bigr|^{p_{\sm}}\biggr)^{\f{1}{p_{\sm}}}<\iy.$$
For $r>0$, we consider the function space of all measurable
functions $f$ on $Q_r$ such that
$$\|f\|_{L^{\iy}_t L^2_x(Q_r)}
:=\sup_{t\in(-r^{\sm},0]}\biggl(\int_{B_r}\bigl|f(x,t)\bigr|^2\biggr)^{\f{1}{2}}<\iy.$$

\begin{thm} Suppose that a function $u\in L^{\iy}_T(L^1_{\om})\cap C(\pa_p Q_{1+\e})$
is a viscosity solution of the equation $$\rL_0 u-\pa_t
u=h\,\,\text{ in $Q_{1+\e}$ }$$ for $h\in L^{\iy}_T(L^2_x)$,
where $\rL_0\in\fL_0(\sm)$ for $\sm\in[\sm_0,2)$ with
$\sm_0\in(1,2)$. Then there is a solution $v\in
L^{\iy}_T({\dot{H}}^{\f{\sm}{2}}_x)$ of the equation $\rL_0 v-\pa_t
v=h\mathbbm{1}_{Q_{1+\e}}$ in $\BR^n_T$ such that
$$\sup_{\rL\in\fL_0(\sm)}\|\rL u\|_{L^{\iy}_t L^2_x(Q_{1/2})}\lesssim
\sup_{Q_{1+\e}}|u-v|+\|u-v\|_{L^{\iy}_T(L^1_{\om})}
+\|h\|_{L^{\iy}_T(L^2_x)}+\|v\|_{L^{\iy}_T({\dot{H}}^{\f{\sm}{2}}_x)}.$$
\end{thm}

\pf Take any $\rL\in\fL_0(\sm)$ for $\sm\in[\sm_0,2)$ with
$\sm_0\in(1,2)$. Let $v\in L^{\iy}_T({\dot{H}}^{\f{\sm}{2}}_x)$ be a
solution of the equation $\rL_0 v-\pa_t v=h\mathbbm{1}_{Q_{1+\e}}$
in $\BR^n_T$. By the Sobolev embedding theorem, we see that
\begin{equation}v\in
L^{\iy}_T({\dot{H}}^{\f{\sm}{2}}_x)\subset
L^{\iy}_T(L^{p_{\sm}}_x).\end{equation} Since $v\in
L^{\iy}_T({\dot{H}}^{\f{\sm}{2}}_x)$ is equivalent to
$(-\Delta)^{\sm/2}v\in L^{\iy}_T(L^2_x)$, it follows from Lemma 4.2
that $\rL_0 v\in L^{\iy}_T(L^2_x)$, and so $\,\pa_t v\in
L^{\iy}_T(L^2_x)$. By H\"older's inequality and (4.3), we have that
$v\in L^{\iy}_T(L^1_{\om})$. From Theorem 4.1, we obtain that
\begin{equation}\|u-v\|_{C^{2,\ap}(Q_1)}\lesssim
\sup_{Q_{1+\e}}|u-v|+\|u-v\|_{L^{\iy}_T(L^1_{\om})}.\end{equation}
Since $\mu_t(u-v,x,y)=\int_0^1\int_0^1\la D^2(u-v)((x+\tau y)-2s\tau
y,t)y,y\ra\,ds\,d\tau$ by the mean value theorem, we have that
$$\bigl|\mu_t(u-v,x,y)\bigr|\lesssim
\bigl(\,\sup_{Q_{1+\e}}|u-v|+\|u-v\|_{L^{\iy}_T(L^1_{\om})}\,\bigr)|y|^2$$
for any $(x,t)\in Q_{1/2}$ and $y\in B_{\f{1}{2}+\e}$. So we get that
\begin{equation}\begin{split}\bigl|\rL(u-v)(x,t)\bigr|&\lesssim
\bigl(\,\sup_{Q_{1+\e}}|u-v|+\|u-v\|_{L^{\iy}_T(L^1_{\om})}\bigr)
\int_{|y|<\f{1}{2}+\e}|y|^2 K(y)\,dy \\&\quad+(u-v)* K_{\e}(y)\lesssim
\sup_{Q_{1+\e}}|u-v|+\|u-v\|_{L^{\iy}_T(L^1_{\om})}
\end{split}\end{equation} for any $(x,t)\in Q_{1/2}$, where
$K_{\e}(y)=\mathbbm{1}_{\BR^n\s B_{\f{1}{2}+\e}}(y)K(y)$. This
implies that $$\|\rL(u-v)\|_{L^{\iy}_t L^2_x(Q_{1/2})}\lesssim
\sup_{Q_{1+\e}}|u-v|+\|u-v\|_{L^{\iy}_T(L^1_{\om})}.$$
Hence we conclude that
$$\sup_{\rL\in\fL_0(\sm)}\|\rL u\|_{L^{\iy}_t L^2_x(Q_{1/2})}\lesssim
\sup_{Q_{1+\e}}|u-v|+\|u-v\|_{L^{\iy}_T(L^1_{\om})}
+\|h\|_{L^{\iy}_T(L^2_x)}+\|v\|_{L^{\iy}_T({\dot{H}}^{\f{\sm}{2}}_x)}.\quad\qed$$

\section{Local uniform upper boundedness of viscosity subsolutions }

In this section, local uniform upper boundedness of viscosity
subsolutions in $L^{\iy}_T(L^1_{\om})$ will be achieved by using
almost the same idea of the proof of the Harnack inequality in
\cite{KL3}.

\begin{thm} Let $\sm\in(1,2)$. If $u\in L^{\iy}_T(L^1_{\om})\cap
C(Q_2)$ satisfies the equation $$\bM_0^+ u-\pa_t
u\ge-\|u\|_{L^{\iy}_T(L^1_{\om})}\,\,\text{ in $Q_2$ }$$ in the
viscosity sense, then we have the estimate
$$\sup_{Q_{1/2}} u\lesssim\|u\|_{L^{\iy}_T(L^1_{\om})}.$$\end{thm}

\pf Without loss of generality, we may assume that
$u\in\rB(\BR^n_T)$. Indeed, if we set $u_1=u\mathbbm{1}_{Q_2}$ and
$u_2=u\mathbbm{1}_{\BR^n_T\s Q_2}$, then it easily follows that
$$\bM_0^+ u_1-\pa_t u_1\gtrsim-\|u\|_{L^{\iy}_T(L^1_{\om})}\,\,\text{ in
$Q_2$. }$$ Since $u$ is continuous on $Q_2$, $u_1$ is bounded on
$\BR^n_T$. So we could use $u_1$ instead of $u$. Also we may assume
that $\|u\|_{L^{\iy}_T(L^1_{\om})}=1$ by dividing $u$ by the norm
$\|u\|_{L^{\iy}_T(L^1_{\om})}$. Thus it suffices to show that
$\,\sup_{Q_{1/2}} u\le C.$ If $u$ is non-positive on $Q_{1/2}$, then
there is nothing to prove it. Thus we may now suppose that $u$ is
non-negative on $Q_{1/2}$. We set $s_0=\inf\{s>0:u(x,t)\le s\,
\dd((x,t),\pa_p Q_1)^{-n-\sm},\,\forall\,(x,t)\in Q_1\}.$ Then we
see that $s_0>0$ and there is some $(\check x,\check t)\in Q_1$ such
that $$u(\check x,\check t)=s_0\,\dd((\check x,\check t),\pa_p
Q_1))^{-n-\sm}=s_0\dd_0^{-n-\sm}$$ where $\dd_0=\dd((\check x,\check
t),\pa_p Q_1)\le 2^{1/\sm}<2$ for $\sm\in(1,2)$. We note that
\begin{equation}\rB^{\dd}_r(x_0,t_0)\subset
Q_r(x_0,t_0)\subset\rB^{\dd}_{2r}(x_0,t_0)\end{equation} for any
$r>0$ and $(x_0,t_0)\in\BR^n_T$.

To finish the proof, we have only to show that $s_0$ can not be too
large because $u(x,t)\le C_1 \dd((x,t),\pa_p Q_1)^{-n-\sm}\le C$ for
any $(x,t)\in Q_{1/2}\subset Q_1$ if $C_1>0$ is some constant with
$s_0\le C_1$. Assume that $s_0$ is very large. Then by Chebyshev's
inequality we have that
$$\bigl|\{u\ge u(\check x,\check t)/2\}\cap
Q_1\}\bigr|\le\f{2}{|u(\check x,\check
t)|}\|u\|_{L^{\iy}(L^1_{\om})}\lesssim s_0^{-1}\dd_0^{n+\sm}.$$ Since
$\rB^{\dd}_r(\check x,\check t)\subset Q_1$ and
$|\rB^{\dd}_{r}|=C\dd_0^{n+\sigma}$ for $r=\dd_0/2\le
2^{-(1-1/\sm)}<1$ for $\sm\in(1,2)$, we easily obtain that
\begin{equation}\bigl|\{u\ge u(\check x,\check t)/2\}\cap
\rB^{\dd}_r(\check x,\check t)\}\bigr|\lesssim
s_0^{-1}|\rB^{\dd}_r|.\end{equation} In order to get a
contradiction, we estimate $|\{u\le u(\check x,\check t)/2\}\cap
\rB^{\dd}_{\dt r/2}(\check x,\check t)|$ for some very small $\dt>0$
(to be determined later). For any $(x,t)\in\rB^{\dd}_{2\dt r}(\check
x,\check t)$, we have that $u(x,t)\le
s_0(\dd_0-\dt\dd_0)^{-n-\sm}=u(\check x,\check t)(1-\dt)^{-n-\sm}$
for $\dt>0$ so that $(1-\dt)^{-n-\sm}$ is close to $1$. We consider
the function
$$v(x,t)=\f{u(\check x,\check t)}{(1-\dt)^{n+\sm}}-u(x,t).$$ Then we see that $v\ge 0$ on
$\rB^{\dd}_{2\dt r}(\check x,\check t)$, and also $\bM_0^- v-\pa_t
v\le 1$ on $Q_{\dt r}(\check x,\check t)$ because $\bM_0^+ u-\pa_t
u\ge -1$ on $Q_{\dt r}(\check x,\check t)$. In order to apply
Theorem 4.12 \cite{KL3} to $v$, we consider $w=v^+$ instead of $v$.
Since $w=v+v^-$, we have that
\begin{equation}\bM_0^- w-\pa_t w\le\bM_0^- v-\pa_t v+\bM_0^+ v^--\pa_t v^-\le
1+\bM_0^+ v^--\pa_t v^-\end{equation} on $Q_{\dt r}(\check x,\check
t).$ Since $v^-\equiv 0$ on $\rB^{\dd}_{2\dt r}(\check x,\check t)$,
if $(x,t)\in Q_{\dt r}(\check x,\check t)$ then we have that
$\mu_t(v^-,x,y)=v^-(x+y,t)+v^-(x-y,t)$ for $y\in\BR^n$.

Take any $(x,t)\in Q_{\dt r}(\check x,\check t)$ and any
$\vp\in\rC^2_{Q_{\dt r}(\check x,\check t)}(v^-;x,t)^+$. Since
$(x,t)+Q_{\dt r}\subset Q_{2\dt r}(\check x,\check t)$ and
$v^-(x,t)=0$, we see that $\pa_t\vp(x,t)=0$. Thus we have that
\begin{equation*}\begin{split}\bM_0^+
v^-(x,t)-\pa_t\vp(x,t)&=(2-\sm)\int_{\BR^n}\f{\Ld\mu_t^+(v^-,x,y)-\ld\mu_t^-(v^-,x,y)}{|y|^{n+\sm}}\,dy\\
&\le 2(2-\sm)\Ld\int_{\{y\in\BR^n:\,v(x+y,t)<0\}}\f{-v(x+y,t)}{|y|^{n+\sm}}\,dy\\
&\le 2(2-\sm)\Ld\int_{B^c_{\dt
r}}\f{\bigl(u(x+y,t)-(1-\dt)^{-n-\sm}u(\check x,\check t)\bigr)_+}{|y|^{n+\sm}}\,dy\\
&\leq  C(2-\sm)\Ld\bigl((\dt
r)^{-n-\sigma}+1\bigr)\int_{\BR^n}\f{|u(y,t)|}{1+|y|^{n+\sm}}\,dy.
\end{split}\end{equation*}
This implies that $$\bM_0^+ v^--\pa_t v^-\lesssim
\|u\|_{L^{\iy}_T(L^1_{\om})}(\dt r)^{-n-\sigma}\lesssim(\dt
r)^{-n-\sigma}\text{ on $Q_{\dt r}(\check x,\check t)$. }$$ Thus by
(5.3), we obtain that $w$ satisfies
$$\bM_0^- w(x,t)-\pa_t w\lesssim (\dt
r)^{-n-\sm}\,\,\text{ on $Q_{\dt r}(\check x,\check t)$ }$$ in
viscosity sense. Since $u(\check x,\check
t)=s_0\dd_0^{-\bt}=2^{-\bt}s_0 r^{-\bt}$, by Theorem 4.12 \cite{KL3}
there is some $\vep_*>0$ such that
\begin{equation*}\begin{split}&\bigl|\{u\le u(\check
x,\check t)/2\}\cap\rB^{\dd}_{\dt r/2}(\check x,\check
t)\bigr|\le\bigl|\{u\le u(\check x,\check t)/2\}\cap Q_{\dt
r/2}(\check x,\check t)\bigr|\\
&\qquad\qquad\qquad=\bigl|\{w\ge u(\check x,\check
t)((1-\dt)^{-\bt}-1/2)\}
\cap Q_{\dt r/2}(\check x,\check t)\bigr|\\
&\qquad\qquad\qquad\lesssim(\dt
r)^{n+\sigma}\bigl[((1-\dt)^{-\bt}-1)u(\check x,\check t)+C(\dt
r)^{-\sm}(\dt
r)^{\sm}\bigr]^{\vep_*}\\
&\qquad\qquad\qquad\qquad\qquad\qquad\qquad\qquad\times\bigl[u(\check x,\check t)((1-\dt)^{-\bt}-1/2)\bigr]^{-\vep_*}\\
&\qquad\qquad\qquad\lesssim(\dt
r)^{n+\sigma}\biggl[\biggl(\f{(1-\dt)^{-\bt}-1}{(1-\dt)^{-\bt}-1/2}\biggr)^{\vep_*}
+\f{s_0^{-\vep_*}r^{n+\sm}}{((1-\dt)^{-\bt}-1/2)^{\vep_*}}\biggr]\\
&\qquad\qquad\qquad\lesssim(\dt
r)^{n+\sigma}[((1-\dt)^{-\bt}-1)^{\vep_*}+s_0^{-\vep_*}r^{n+\sm}].
\end{split}\end{equation*}
We now choose $\dt>0$ so small enough that $(\dt
r)^{n+\sigma}((1-\dt)^{-\bt}-1)^{\vep_*}\lesssim |\rB^{\dd}_{\dt
r/2}|/4.$ Since $\dt$ was chosen independently of $s_0$, if $s_0$ is
large enough for such fixed $\dt$ then we get that $(\dt
r)^{n+\sm}s_0^{-\vep_*}r^{n+\sm}\lesssim |\rB^{\dd}_{\dt r/2}|/4.$
Therefore we obtain that $$\bigl|\{u\le u(\check x,\check
t)/2\}\cap\rB^{\dd}_{\dt r/2}(\check x,\check t)\bigr|\le
|\rB^{\dd}_{\dt r/2}|/2.$$ Thus we conclude that
\begin{equation*}\begin{split}\bigl|\{u\ge u(\check x,\check t)/2\}\cap\rB^{\dd}_r(\check x,\check t)\bigr|
&\ge\bigl|\{u\ge u(\check x,\check t)/2\}\cap \rB^{\dd}_{\dt
r/2}(\check x,\check t)\bigr|\\&\ge\bigl|\{u>u(\check x,\check
t)/2\}\cap\rB^{\dd}_{\dt r/2}(\check x,\check
t)\bigr|\\&\ge\bigl|\rB^{\dd}_{\dt r/2}(\check x,\check
t)\bigr|-\bigl|\rB^{\dd}_{\dt r/2}\bigr|/2\\&=\bigl|\rB^{\dd}_{\dt
r/2}\bigr|/2=C |\rB^{\dd}_r|,\end{split}\end{equation*} which
contradicts (5.2) if $s_0$ is large enough. Hence we complete the
proof. \qed

\section{Local uniform boundedness of linear operators}

The main theme of this section is to establish local uniform
boundedness of linear operators from the result obtained in Section
5, which facilitate obtaining local uniform boundedness of extremal
operators to be given in the next section.

\begin{lemma} Let $u\in L^{\iy}_T(L^1_{\om})\cap C(\pa_p Q_2)$ be a
viscosity solution satisfying the equation
$$\bI u-\pa_t u=0\,\text{ in
$Q_2$.}$$ If $K\in\cK_0(\sm)$ is a symmetric kernel,
then for any cut-off
function $\vp\in C^{\iy}_c(\BR^n)$ supported in $B_1$ and with
$0\le\vp\le 1$ in $\BR^n,$ we have that $$\bM^+_2 u_{\vp}-\pa_t
u_{\vp}\ge 0\,\,\text{ in $Q_1$ }$$ in the viscosity sense, where
$$u_{\vp}(x,t)=\int_{\BR^n}\mu_t(u,x,y)K(y)\vp(y)\,dy.$$
\end{lemma}

\pf By Lemma 5.4 \cite{KL4} and Lemma 3.2, without loss of
generality we may assume that $u\in C^2(Q_1)$. So we see that
integro-differential type operators like
$u_{\vp}$ are well-defined and continuous in $Q_1$. For $\el\in\BN$,
we set $\vp_{\el}(y)=\mathbbm{1}_{\BR^n\s B_{4/\el}}(y)K(y)\vp(y)$.
Then we see that $\vp_{\el}\in L^1(\BR^n)$ for all $\el\in\BN$. By
Lebesgue's dominated convergence theorem, we have that
$$u_{\vp}=\lim_{\el\to\iy}\int_{\BR^n}\mu_\cdot(u,\cdot,y)\vp_{\el}(y)\,dy
=2\bigl(\lim_{\el\to\iy}u*\vp_{\el}-u\,\|\vp_{\el}\|_{L^1}\bigr).$$
Now it follows from Lemma 3.4 that
$$\bI\biggl(u*\f{\vp_{\el}}{\|\vp_{\el}\|_{L^1}}\biggr)
-\pa_t\biggl(u*\f{\vp_{\el}}{\|\vp_{\el}\|_{L^1}}\biggr)\ge
0\,\,\text{ in $Q_1$. }$$ Also we have that $\bI u-\pa_t u=0$ in
$Q_1$. Thus by applying Theorem 2.4 \cite{KL3}, we easily obtain
that
\begin{equation*}\begin{split}&\bM^+_2\bigl(u*\vp_{\el}-u\,\|\vp_{\el}\|_{L^1}\bigr)
-\pa_t\bigl(u*\vp_{\el}-u\,\|\vp_{\el}\|_{L^1}\bigr)\\
&=\|\vp_{\el}\|_{L^1}\biggl[\bM^+_2\biggl(u*\f{\vp_{\el}}{\|\vp_{\el}\|_{L^1}}-u\biggr)
-\pa_t\biggl(u*\f{\vp_{\el}}{\|\vp_{\el}\|_{L^1}}-u\biggr)\biggr]\ge
0\,\text{ in $Q_1$ }
\end{split}\end{equation*} for any $\el\in\BN$. Hence
we can obtain the required result by taking limit $\el\to\iy$. \qed

\begin{lemma} Let $u\in L^{\iy}_T(L^1_{\om})\cap C(\pa_p Q_2)$ be any
viscosity solution satisfying the equation
$$\bI u-\pa_t u=0\,\text{ in $Q_2$.}$$
Then we have the estimate
$$\bM^+_2(\rL u)-\pa_t(\rL
u)\gtrsim-\|u\|_{L^{\iy}_T(L^1_{\om})}\,\,\text{ in $Q_1$ }$$ for any
$\rL\in\fL_2$.
\end{lemma}

\pf By Lemma 5.4 \cite{KL4} and Lemma 3.2, without loss of
generality we may assume that $u\in C^2(Q_1)$. For $\el\in\BN$, let $\e_{\el}(y)=\mathbbm{1}_{\BR^n\s
B_{4/\el}}(y)K(y)$. Take any $\rL\in\fL_2$. Then as in Lemma 6.1 we
have that
$$\rL u=\lim_{\el\to\iy}\int_{\BR^n}\mu_\cdot(u,\cdot,y)\e_{\el}(y)\,dy
=2\lim_{\el\to\iy}\bigl(u*{\e_{\el}}-u\,\|\e_{\el}\|_{L^1}\bigr).$$
Let $\vp\in C^{\iy}_c(\BR^n)$ be any radial cut-off function
supported in $B_2$ such that $\vp\equiv 1$ in $B_{3/2}$ and
$0\le\vp\le 1$ in $\BR^n$. We set $\phi_{\el}(y)=\e_{\el}(y)\vp(y)$
and $\psi_{\el}(y)=\e_{\el}(y)(1-\vp(y))$. By Lemma 6.1, we have
that
\begin{equation}\bM^+_2\bigl(u*\phi_{\el}-u\,\|\phi_{\el}\|_{L^1}\bigr)
-\pa_t\bigl(u*\phi_{\el}-u\,\|\phi_{\el}\|_{L^1}\bigr)\ge
0\,\,\text{ in $Q_1$. }
\end{equation}

Also we now estimate $\bI(u*\psi_{\el})-\pa_t(u*\psi_{\el})$ in
$Q_1$. Take any point $(x,t)\in Q_1$. We note that
\begin{equation*}\begin{split}\bI(u*\psi_{\el})-\pa_t(u*\psi_{\el})
&=\inf_{\bt}\rL_{\bt}(u*\psi_{\el})-\pa_t(u*\psi_{\el})\\
&=\inf_{\bt}u*(\rL_{\bt}\psi_{\el})-\pa_t(u*\psi_{\el})
\end{split}\end{equation*}
and
\begin{equation*}\begin{split}u*\rL_{\bt}(\psi_{\el})(x,t)&=\int_{\BR^n}u(x-y,t)\int_{|z|\ge
\f{1}{2}}\mu(\psi_{\el},y,z)K(z)\,dz\,dy\\
&\quad+\int_{|y|\ge 1}u(x-y,t)\int_{|z|<
\f{1}{2}}\mu(\psi_{\el},y,z)K(z)\,dz\,dy\\
&:=I(x,t)+II(x,t)
\end{split}\end{equation*} by the definition of $\psi_{\el}$. Then
it is easy to check that
\begin{equation}I=2(u*\psi_{\el}*\e_2)-2 c\,(u*\psi_{\el})
\end{equation} for a universal constant $c>0$. By the mean value theorem and
triangle inequality, we see that for any $y\in\BR^n\s B_1$ and $z\in
B_{1/2}$,
$$\mu(\psi_{\el},y,z)=\int_0^1\int_0^1\la D^2\psi_{\el}((y+\tau
z)-2s\tau z) z,z\ra\,ds\,d\tau,$$
$$|(y+\tau z)-2s\tau z|=|y+\tau(1-2s)z|\ge|y|-|z|\ge|y|/2.$$
Since
$D^2\psi_{\el}=(D^2\e_{\el})(1-\vp)-2(D\e_{\el})(D\vp)-\e_{\el}(D^2\vp)$,
by (1.2) and (1.4) we have that $$|D^2\psi_{\el}((y+\tau z)-2s\tau
z,t)|\le\f{C}{|y|^{n+\sm}}\mathbbm{1}_{\BR^n\s B_3}(y):=k(y)$$ for
any $y\in\BR^n\s B_1$, $z\in B_{1/2}$ and $s,\tau\in[0,1]$. Thus we
obtain that
\begin{equation}|II(x,t)|\le
|u|*k(x,t)\int_{|z|<\f{1}{2}}|z|^2 K(z)\,dz\lesssim|u|*k(x,t)\lesssim
\|u\|_{L^{\iy}_T(L^1_{\om})}.
\end{equation} Hence it easily
follows from (6.2), (6.3) and Young's inequality that
\begin{equation}|u*\rL_{\bt}(\psi_{\el})(x,t)|\lesssim\|u\|_{L^{\iy}_T(L^1_{\om})}
\end{equation}
for any $\bt$, and thus we have that
\begin{equation}\rI(u*\psi_{\el})(x,t)\gtrsim -\|u\|_{L^{\iy}_T(L^1_{\om})}.
\end{equation}
Since $u\in C^1_t(Q_1)$, as in the above estimate we can obtain that
\begin{equation}\begin{split}\pa_t(u*\psi_{\el})(x,t)&=(\pa_t
u)*\psi_{\el}(x,t)=(\rI u)*\psi_{\el}(x,t)\\
&\le(\rL_{\bt}u)*\psi_{\el}(x,t)=u*(\rL_{\bt}\psi_{\el})(x,t)\\
&\lesssim\|u\|_{L^{\iy}_T(L^1_{\om})}.
\end{split}\end{equation} Hence by (6.1), (6.5) and (6.6), we
conclude that
$$\bM^+_2(\rL u)-\pa_t(\rL
u)\gtrsim-\|u\|_{L^{\iy}_T(L^1_{\om})}\,\,\text{ in $Q_1$. }$$
Therefore we complete the proof. \qed

\begin{lemma} If $u\in L^{\iy}_T(L^1_{\om})\cap C(\pa_p Q_{2+\e})$ is a
viscosity solution satisfying the equation
$$\bI u-\pa_t u=0\,\text{ in
$Q_{2+\e}$}$$ where $\bI$ is defined on $\fL_2(\sm)$ for $\sm\in(\sm_0,2)$ with $\sm_0\in(1,2)$, then we have the estimate
$$\sup_{\rL\in\fL_2}\biggl(\sup_{Q_{1/4}}\rL u\biggr)\lesssim\|u\|_{L^{\iy}_T(L^1_{\om})}.$$
\end{lemma}

\pf By Lemma 5.4 \cite{KL4} and Lemma 3.2, without loss of
generality we may assume that $u\in C^2(Q_2)$. By
Lemma 6.2, we see that
\begin{equation}\bM^+_0(\rL u)-\pa_t(\rL u)\ge\bM^+_2(\rL
u)-\pa_t(\rL u)\gtrsim-\|u\|_{L^{\iy}_T(L^1_{\om})}\,\,\text{ in $Q_1$
}\end{equation} for any $\rL\in\fL_2$. Since it is easy to check
that $\rL$ is a nonlocal parabolic operator, we see that $\rL u\in
C(Q_2)$ (see \cite{KL4}).

Let $\vp\in C^{\iy}_c(\BR^n)$ be a function such that $0\le\vp\le
1$, $\vp=1$ in $B_2$, $\vp=0$ in $\BR^n\s B_{2+\e/2}$ and
$|D^2\vp|\le N_0$ in $B_{2+\e/2}$ for some $N_0>0$. Then by the
change of variables we have that
\begin{equation}\int_{\BR^n}\rL
u(x,t)\,\vp(x)\,dx=\int_{\BR^n}u(x,t)\,\rL\vp(x)\,dx.
\end{equation} We note that $|(x+\tau y)-2s\tau
y|=|x+\tau(1-2s)y|\le|x|+|y|$ for $s,\tau\in[0,1]$ and
$$\mu(\vp,x,y)=\int_0^1\int_0^1\la D^2\vp((x+\tau
y)-2s\tau y) y,y\ra\,ds\,d\tau$$ for any $x\in B_1$ and $y\in B_1$.
We now have that
\begin{equation*}\rL\vp(x)=\int_{B_1}\mu(\vp,x,y)K(y)\,dy+\int_{\BR^n\s B_1}\mu(\vp,x,y)K(y)\,dy
:=b(x)+c(x)\end{equation*} and $c(x)=2\,\vp*\e_4(x)-2\,c_0\vp(x)$
where $c_0=\int_{\BR^n\s B_1}K(y)\,dy<\iy$ and
$\e_r(y)=\mathbbm{1}_{\BR^n\s B_r}(y)K(y)$. Then it is easy to check
that $|b(x)|\le N_0\int_{B_1}|y|^2 K(y)\,dy\le c<\iy$ for $|x|<5$
and $|b(x)|=0$ for $|x|\ge 5$, and $|c(x)|\le c$ for $|x|<5$ and
$|c(x)|\le c/|x|^{n+\sm}$ for $|x|\ge 5$, where $c>0$ is a universal
constant. So we see that $|\rL\vp(x)|\lesssim\om(x)$. Thus by (6.8), we obtain that
\begin{equation}\biggl|\int_{\BR^n}\rL
u(x,t)\,\vp(x)\,dx\biggr|\lesssim
\|u\|_{L^{\iy}_T(L^1_{\om})}.\end{equation} We set $\phi(x)=1-\vp(x)$ and $w(x,t)=\vp(x)\,\rL
u(x,t)$, and we denote by $f^x(y)=f(x+y)$. Then (6.9) implies that
$w\in L^{\iy}_T(L^1_{\om})\cap C(Q_2)$. We now estimate $\bM_0^+
w(x,t)$ for $x\in B_1$ and $t\in(-T,0]$. For this, as in (6.4) we
have that
\begin{equation}\sup_{(x,t)\in Q_1}\bigl|u*\rL(\phi^x
K)(x,t)\bigr|\lesssim\|u\|_{L^{\iy}_T(L^1_{\om})},
\end{equation} because $\phi^x K$ is a smooth function with nice decay such that $\phi^x K=0$
on $B(x;1)$ for each $x\in B_1$. If $(x,t)\in Q_1$, then by the
change of variables and (6.10), we have the estimate
\begin{equation}\begin{split}\rL_{\bt}w(x,t)
&=\int_{\BR^n}\mu_t(\rL u,x,y)K(y)\,dy-\int_{\BR^n}\mu_t\bigl((\rL
u)\phi,x,y\bigr)K(y)\,dy\\
&=\int_{\BR^n}\mu_t(\rL u,x,y)K(y)\,dy-2\int_{\BR^n}\rL
u(x+y,t)\phi^x(y) K(y)\,dy\\
&=\int_{\BR^n}\mu_t(\rL
u,x,y)K(y)\,dy-2\,u*\rL(\phi^x K)(x,t)\\
&\gtrsim\int_{\BR^n}\mu_t(\rL
u,x,y)K(y)\,dy-\,\|u\|_{L^{\iy}_T(L^1_{\om})}
\end{split}\end{equation}
for any $\rL_{\bt}\in\fL_2$. Hence by (6.7) and (6.11) we conclude
that
\begin{equation*}\bM_0^+ w-\pa_t w\gtrsim\bM^+_2(\rL
u)-\pa_t(\rL
u)-\,\|u\|_{L^{\iy}_T(L^1_{\om})}\gtrsim-\|u\|_{L^{\iy}_T(L^1_{\om})}\,\,\text{
on $Q_1$. }
\end{equation*}
Therefore the required result can be achieved by applying Theorem
5.1. \qed

\section{Local uniform boundedness of extremal operators}

In this section, we show that if $u\in L^{\iy}_T(L^1_{\om})$ is a
viscosity solution of the nonlocal parabolic concave equation $\bI
u-\pa_t u=0$ in $Q_2$, then $\bM^+_0 u$ and $\bM^-_0 u$ are bounded
uniformly on $Q_{1/2}$ for $\sm\in(\sm_0,2)$ with $\sm_0\in(1,2)$.
This plays an important role as a cornerstone in proving the main
theorem in the final section.

\begin{lemma} Let $u\in L^{\iy}_T(L^1_{\om})\cap C(\pa_p Q_2)$ be a
viscosity solution satisfying the equation
$$\bI u-\pa_t u=0\,\text{ in
$Q_2$.}$$ If $K$ is a symmetric kernel with
$K(y)\le(2-\sm)\Ld|y|^{-n-\sm}$, then for any function $\gm\in
C^{\iy}_c(-T,T]$ with $\gm=1$ in $(-2^{1-\sm},0]$ and
$\supp(\gm)\subset(-1,\e]$ and any radial cut-off function $\psi\in
C^{\iy}_c(\BR^n)$ supported in $B_2$ such that $\psi=1$ in
$B_{8/5}$, $\psi=0$ in $\BR^n\s B_2$ and $0\le\psi\le 1$ in $\BR^n,$
we have that
$$\bM^+_2 (\psi\gm\, u_{\vp})-\pa_t(\psi\gm\, u_{\vp})\gtrsim-\|u\|_{L^{\iy}_T(L^1_{\om})}
\,\,\text{ in $Q_{1/2}$ }$$ in the viscosity sense, where
$$u_{\vp}(x,t)=\int_{\BR^n}\mu_t(u,x,y)K(y)\vp(y)\,dy$$
for a radial cut-off function $\vp\in C^{\iy}_c(\BR^n)$ supported in
$B_{1/4}$ with $0\le\vp\le 1$ in $\BR^n.$
\end{lemma}

\pf By Lemma 6.1, we see that  $\bM^+_2 u_{\vp}-\pa_t u_{\vp}\ge 0$
in $Q_1$ in the viscosity sense. Set $\phi=1-\psi\gm$ in $\BR^n_T$.
Take any $\rL_{\bt}\in\fL_2$ and $(x,t)\in Q_{1/2}$. Then we have
that
\begin{equation}\begin{split} \rL_{\bt}(\psi\gm\,
u_{\vp})(x,t)=\int_{\BR^n}\mu_t(u_{\vp},x,y)K_{\bt}(y)\,dy-E(x,t)
\end{split}\end{equation} where $E(x,t)=\int_{\BR^n}\mu_t(\phi\,
u_{\vp},x,y)K_{\bt}(y)\,dy$. By the mean value theorem and triangle
inequality, we see that
\begin{equation}\mu(\phi^x
K_{\bt},y,z)=\int_0^1\int_0^1\la D^2(\phi^x K_{\bt})((y+\tau
z)-2s\tau z) z,z\ra\,ds\,d\tau
\end{equation} and $|x+(y+\tau
z)-2s\tau
z|=|x+y+\tau(1-2s)z|\ge|y|-\f{15}{16}\,|y|\ge\f{1}{16}\,|y|$ for any
$y\in\BR^n\s B_{4/5}$, $z\in B_{1/4}$ and $-2^{1-\sm}<t\le 0$. Also
we note that $\mu(\phi^x K_{\bt},y,z)=0$ for any $y\in B_{4/5}$,
$z\in B_{1/4}$ and $-2^{1-\sm}<t\le 0$. Thus by (7.2) we obtain that
\begin{equation*}\begin{split}
E(x,t)&=2\int_{\BR^n}\biggl(\int_{\BR^n}\mu_t(u,x+y,z)K(z)\vp(z)\,dz\biggr)
\phi(x+y) K_{\bt}(y)\,dy\\
&=2\int_{\BR^n}\biggl(\int_{\BR^n}\mu_t(u,x+y,z)\phi(x+y)
K_{\bt}(y)\,dy\biggr)K(z)\vp(z)\,dz\\
&=2\int_{\BR^n}u(x+y,t)\biggl(\int_{\BR^n}\mu_t(\phi^x
K_{\bt},y,z)K(z)\vp(z)\,dz\biggr)dy\\
&=2\int_{\BR^n}u(x+y,t)\biggl(\int_{|z|<\f{1}{4}}\mu_t(\phi^x
K_{\bt},y,z)K(z)\vp(z)\,dz\biggr)dy\\
&\lesssim
\int_{|y|\ge\f{4}{5}}|u(x+y,t)|\f{1}{|y|^{n+2+\sm}}\,dy\int_{|z|<\f{1}{4}}|z|^2
K(z)\,dz\lesssim\|u\|_{L^{\iy}_T(L^1_{\om})}
\end{split}\end{equation*} for any $|x|<1/2$ and $-2^{1-\sm}<t\le 0$.
Hence by (7.1) we conclude that
\begin{equation*}\begin{split}&\bM^+_2 (\psi\gm\,
u_{\vp})(x,t)-\pa_t(\psi\gm\, u_{\vp})(x,t)\\&\qquad\ge\bM^+_2
u_{\vp}(x,t)-\pa_t u_{\vp}(x,t)-E(x,t)\gtrsim -\|u\|_{L^{\iy}_T(\om)}
\end{split}\end{equation*} for any $(x,t)\in Q_{1/2}$. Therefore we
complete the proof. \qed

\begin{lemma} Let $u\in L^{\iy}_T(L^1_\om)\cap C(\pa_p Q_2)$ be any
viscosity solution satisfying the equation $$\bI u-\pa_t u=0\,\,\text{ in
$Q_2$,}$$ where $\bI$ is defined on $\fL_2(\sm)$ for $\sm\in(\sm_0,2)$
with $\sm_0\in(1,2)$. Then for any operator $\rL$ with a symmetric
kernel $K$ satisfying $K(y)\le(2-\sm)\Ld|y|^{-n-\sm}$, we have the estimate
$$\sup_{Q_{1/2}}|\rL u|\lesssim\|u\|_{L^{\iy}_T(L^1_\om)}.$$
\end{lemma}

\pf Take any $\sm\in(\sm_0,2)$ with $\sm_0\in(1,2)$. As in Lemma
6.3, without loss of generality, we may assume that $u\in
C^2(Q_1)$. For convenience, we normalize
$\|u\|_{L^{\iy}_T(L^1_\om)}=1$. By Lemma 6.3, we see that
$\sup_{\rL_{\bt}\in\fL_2}|\rL_{\bt}u|$ is bounded in $Q_{1/2}$ because $-u$ is another
viscosity solution of our equation. So this implies that $|\pa_t u|=|\bI u|$ is bounded in $Q_{1/2}$. Thus it follows from that
\begin{equation*}\begin{split}\|\rL_{\bt}u-\pa_t u\|_{L^{\iy}_t L^2_x(Q_{1/2})}
&\le\|\rL_{\bt}u\|_{L^{\iy}_t L^2_x(Q_{1/2})}+\|\pa_t u\|_{L^{\iy}_t
L^2_x(Q_{1/2})}<\iy.
\end{split}\end{equation*}
Combining Theorem 4.3 with this yields that
\begin{equation}\sup_{\rL\in\fL_0(\sm)}\|\rL u\|_{L^{\iy}_t
L^2_x(Q_{1/2})}<\iy.
\end{equation}
Take any operator $\rL$ with a symmetric kernel $K$ satisfying
$K(y)\le(2-\sm)\Ld|y|^{-n-\sm}$. Then we split $\rL u$ into two
integrals
\begin{equation*}\begin{split}\rL
u(x,t)&=\int_{\BR^n}\mu_t(u,x,y)K(y)\vp(y)\,dy+\int_{\BR^n}\mu_t(u,x,y)K(y)(1-\vp(y))\,dy\\
&:=u_{\vp}(x,t)+u_{1-\vp}(x,t),
\end{split}\end{equation*} where $\vp\in C_c^{\iy}(\BR^n)$ is a radial
cut-off function supported in $B_1$ such that $\vp=1$ in $B_{1/2}$
and $0\le\vp\le 1$ in $\BR^n$. Since $K\in L^1(\BR^n\s B_{1/2})$, it
is easy to check that $\sup_{Q_{1/2}}|u_{1-\vp}|<\iy$, and thus we
have that $$\|u_{1-\vp}\|_{L^{\iy}_t L^2_x(Q_{1/2})}<\iy.$$ Thus by
(7.3), we obtain that
\begin{equation}\|u_{\vp}\|_{L^{\iy}_t L^2_x(Q_{1/2})}<\iy.
\end{equation}
From Lemma 6.1, we have that
\begin{equation}\bM_2^+ u_{\vp}-\pa_t
u_{\vp}\ge 0\,\,\text{ in $Q_1$. }\end{equation} Let $\psi\in
C_c^{\iy}(\BR^n)$ be a function such that $\psi=1$ in $B_{1/2}$ and
$\supp(\psi)\subset B_{\f{1}{2}+\e}$, and let $\gm\in C_c(-T,T]$ be
a function such that $\gm=1$ in $(-2^{-\sm},0]$ and
$\supp(\gm)\subset(-2^{-\sm}-\e,\e]$. Set
$v_{\vp}(x,t)=\psi(x)\gm(t)\,u_{\vp}(x,t)$. Then by (7.4) it is easy
to check that $v_{\vp}\in L^{\iy}_T(L^1_{\om})$. So it follows from
Lemma 7.1 that
\begin{equation*}\bM_2^+v_{\vp}-\pa_t
v_{\vp}\gtrsim -1\,\text{ in $Q_{1/2}$. }\end{equation*} Applying
Theorem 5.1, we obtain that $v_{\vp}\lesssim 1$ in $Q_{1/8}$. Thus the
required upper bound for $\rL u$ on $Q_{1/2}$ follows from a
standard covering and scaling argument.

For the lower bound for $\rL u$ on $Q_{1/2}$, we take an operator
$\rL_{\bt}\in\fL_2(\sm)$ with kernel $K_{\bt}$ and consider an
operator $\rL_*$ with kernel $K_*=\f{2}{\ld}K_{\bt}-\f{1}{\Ld}K$.
Then it is easy to check that $$\f{2-\sm}{|y|^{n+\sm}}\le
K_*(y)\le\f{(2-\sm)(\f{2\Ld}{\ld}-\f{\ld}{\Ld})}{|y|^{n+\sm}}.$$ As
in the first half, we obtain that $\rL_* u\lesssim 1$ in $Q_{1/2}$. This
implies that $\rL u\gtrsim -1$ in $Q_{1/2}$. Therefore the required
result can be achieved. \qed

From the above result, it is natural to obtain the following
corollaries.

\begin{cor} Let $u\in L^{\iy}_T(L^1_\om)\cap C(\pa_p Q_2)$ be any
viscosity solution satisfying the equation $\bI u-\pa_t u=0$ in
$Q_2$, where $\bI$ is defined on $\fL_2(\sm)$ for $\sm\in(\sm_0,2)$
with $\sm_0\in(1,2)$. Then $\bM^+_0 u$, $\bM^-_0 u$ and $\pa_t u$
are uniformly bounded in $Q_{1/2}$, and moreover we have
$$\bigl(\,\sup_{Q_{1/2}}|\bM^+_0 u|\bigr)\vee
\bigl(\,\sup_{Q_{1/2}}|\bM^-_0 u|\bigr)\vee
\bigl(\,\sup_{Q_{1/2}}|\pa_t u|\bigr)\lesssim\|u\|_{L^{\iy}_T(L^1_\om)}.$$
\end{cor}

\begin{cor} Let $u\in L^{\iy}_T(L^1_\om)\cap C(\pa_p Q_2)$ be any
viscosity solution satisfying the equation $\bI u-\pa_t u=0$ in
$Q_2$, where $\bI$ is defined on $\fL_2(\sm)$ for $\sm\in(\sm_0,2)$
with $\sm_0\in(1,2)$. Then we have that
$$\sup_{Q_{1/2}}\int_{\BR^n}\bigl|\mu_{\cdot}(u,\,\cdot\,,y)\bigr|\f{2-\sm}{|y|^{n+\sm}}\,dy
\lesssim\|u\|_{L^{\iy}_T(L^1_\om)}.$$
\end{cor}

\section{Proof of the Main Theorem}

Let $u\in L^{\iy}_T(L^1_\om)$ be any viscosity solution satisfying
the equation
\begin{equation}\bI u-\pa_t u=0\,\text{ in
$Q_2$,}\end{equation} where $\bI$ is defined on $\fL_2(\sm)$ for
$\sm\in(\sm_0,2)$ with $\sm_0\in(1,2)$. From Corollary 7.4, there is a universal constant $c_0>0$ such that
\begin{equation}
\sup_{Q_{1/2}}\int_{\BR^n}|\mu_{\cdot}(u,\,\cdot\,,y)|\f{2-\sm}{|y|^{n+\sm}}\,\vp(y)\,dy\le c_0\,\|u\|_{L^{\iy}_T(L^1_\om)},
\end{equation}
where $\vp\in C^{\iy}_c(\BR^n)$ is a function such that $\vp=1$ in
$B_1$, $\vp=0$ in $\BR^n\s B_{3/2}$ and $0\le\vp\le 1$ in $\BR^n$.

In order to prove Theorem 1.1, our main goal is to obtain that there
is some $\ap\in(0,1)$ such that
\begin{equation}\int_{\BR^n}|\mu_t(u,x,y)-\mu_0(u,0,y)|\f{2-\sm}{|y|^{n+\sm}}\,\vp(y)\,dy\lesssim(|x|+|t|^{\sm})^{\f{\ap}{\sm}}\,\|u\|_{L^{\iy}_T(L^1_\om)}
\end{equation} for any $(x,t)\in Q_{1/2}$. This implies that the
fractional Laplacian $(-\Delta)^{\sm/2}$ admits the H\"older
continuity, and moreover the viscosity solutions of the nonlocal
parabolic equation in Theorem 1.1 enjoy the
$C^{\sm+\ap}$-regularity.

Let $\psi\in C_c^{\iy}(\BR^n)$ be a function such that $\psi=1$ in
$B_{1/2}$ and $\supp(\psi)\subset B_{\f{1}{2}+\e}$, and let $\gm\in
C_c(-T,T]$ be a function such that $\gm=1$ in $(-2^{-\sm},0]$ and
$\supp(\gm)\subset(-2^{-\sm}-\e,\e]$. Set
$w_{\vp}(x,t)=\psi(x)\gm(t)\,v_{\vp}(x,t)$, where
$$v_{\vp}(x,t)=\int_{\BR^n}\bigl[\mu_t(u,x,y)-\mu_0(u,0,y)\bigr]\,\f{2-\sm}{|y|^{n+\sm}}\,\vp(y)\,dy$$
for a radial cut-off function $\vp\in C^{\iy}_c(\BR^n)$ supported in
$B_{1/4}$ with $0\le\vp\le 1$ in $\BR^n.$ Then, as in Lemma 7.2, it
is easy to check that $w_{\vp}\in L^{\iy}_T(L^1_\om)$ and it follows
from Lemma 7.1 that
$$\bM^+_2 w_{\vp}-\pa_t w_{\vp}\gtrsim -\|u\|_{L^{\iy}_T(L^1_\om)}\,\,\text{ in $Q_{1/2}$,}$$
We set
$$v_{\vp}^{\pm}(x,t)=\int_{\BR^n}\bigl[\mu_t(u,x,y)-\mu_0(u,0,y)\bigr]^{\pm}\,\f{2-\sm}{|y|^{n+\sm}}\,\vp(y)\,dy$$
and set
$$w_{\vp}^S(x,t)=\psi(x)\gm(t)\,\int_{\BR^n}\bigl[\mu_t(u,x,y)-\mu_0(u,0,y)\bigr]\,\f{2-\sm}{|y|^{n+\sm}}
\,\vp(y)\mathbbm{1}_S(y)\,dy$$
for a symmetric set $S\subset\BR^n$ (i.e. $S=-S$). Also we consider
the positive part $\rP u$ and negative part $\rN u$ of $w_{\vp}$
defined by $\rP u(x,t)=\psi(x)\gm(t)\,v_{\vp}^+(x,t)$ and $\rN
u(x,t)=\psi(x)\gm(t)\,v_{\vp}^-(x,t)$. Then we see that $\rP
u=\sup_S w_{\vp}^S$ and $\rN u=-\inf_S w_{\vp}^S$, and moreover $\rP
u=w_{\vp}^{S_0}$ and $\rN u=-w_{\vp}^{S_0^c}$ where $S_0$ is the
symmetric set given by
$S_0=\{y\in\BR^n:\mu_t(u,x,y)>\mu_0(u,0,y)\}$.

\begin{lemma} If $u\in L^{\iy}_T(L^1_\om)\cap C(\pa_p Q_2)$ be a viscosity solution
of the equation
\begin{equation*}\bI u-\pa_t u=0\,\text{ in
$Q_2$,}\end{equation*} where $\bI$ is defined on $\fL_2(\sm)$ for
$\sm\in(\sm_0,2)$ with $\sm_0\in(1,2)$, then there exists some $\ap\in(0,1)$ such that
$$\sup_{(x,t)\in Q_{1/8}}\f{\rP u(x,t)}{(|x|^{\sm}+|t|)^{\f{\ap}{\sm}}}\lesssim\|u\|_{L^{\iy}_T(L^1_\om)}.$$
\end{lemma}

\pf We may assume that
$\|u\|_{L^{\iy}_T(L^1_\om)}\le 1$ by dividing the equation by
$\|u\|_{L^{\iy}_T(L^1_\om)}$. Take
any $(x,t)\in Q_{1/8}$ and $\rL\in\fL_2(\sm)$. Then we have that
\begin{equation}\begin{split}\rL(\btau^t_x u-u)(0,0)&=\int_{\BR^n}\bigl[\mu_t(u,x,y)-\mu_0(u,0,y)\bigr]\vp(y)K(y)\,dy\\
&+\int_{\BR^n}\bigl[\mu_t(u,x,y)-\mu_0(u,0,y)\bigr]\phi(y)K(y)\,dy\\
&:=\rL_{\vp}u(x,t)+\rL_{\phi}u(x,t),
\end{split}\end{equation} where $\phi=1-\vp$. Then we see that
\begin{equation}w^-_{\vp}(x,t)\le \rL_{\vp}(x,t)\le w^+_{\vp}(x,t)
\end{equation}
where $w^-_{\vp}(x,t)=\ld\,\rP u(x,t)-\Ld\,\rN u(x,t)$ and
$w^+_{\vp}(x,t)=\Ld\,\rP u(x,t)-\ld\,\rN u(x,t)$. By easy
calculation, the second term in the right hand side of (8.4) becomes
\begin{equation*}\begin{split}\rL_{\phi}u(x,t)
&=2\int_{\BR^n}u(y,t)\bigl[K(y-x)\phi(y-x)-K(y)\phi(y)\bigr]\,dy\\
&\qquad+C^{\phi}_K(t)+2\bigl[u(0,0)-u(x,t)\bigr]\int_{\BR^n}K(y)\phi(y)\,dy
\end{split}\end{equation*} where $C^{\phi}_K(t)=\ds 2\int_{\BR^n}[u(y,t)-u(y,0)]K(y)\phi(y)\,dy$.
Thus it follows from (1.3) and Theorem 3.4 [KL4] that
\begin{equation}\begin{split} A(x,t)&\le\inf_{K\in\cK_2}\bigl[\rL(\btau^t_x u-u)(0,0)-C^{\phi}_K(t)\bigr]\\
&\le\sup_{K\in\cK_2}\bigl[\rL(\btau^t_x u-u)(0,0)-C^{\phi}_K(t)\bigr]\le B(x,t)
\end{split}\end{equation}
for some universal constants $c,\bt>0$, where $A(x,t)=w_{\vp}^-(x,t)-c\,(|x|^{\sm}+|t|)^{\bt/\sm}$ and $B(x,t)=w_{\vp}^+(x,t)+c\,(|x|^{\sm}+|t|)^{\bt/\sm}$. Here we note that $\bt$ could
be chosen freely in the open interval $(0,1)$ (see \cite{KL3}).
Then we have
only three possible cases; either (a) $A(x,t)\le 0$ and $B(x,t)\ge
0$, or (b) $A(x,t)\ge 0$ and $B(x,t)\ge 0$, or (c) $A(x,t)\le 0$ and
$B(x,t)\le 0$.

({\bf  Case I :} (a) $A(x,t)\le 0$ and $B(x,t)\ge 0$ )  (a) implies
that
\begin{equation}\f{\ld}{\Ld}\,\rN
u(x,t)-c_1\,(|x|^{\sm}+|t|)^{\f{\bt}{\sm}}\le\rP
u(x,t)\le\f{\Ld}{\ld}\,\rN
u(x,t)+c_1\,(|x|^{\sm}+|t|)^{\f{\bt}{\sm}}
\end{equation} for any $(x,t)\in
Q_{1/8}$, where $c_1=c/\Ld$.

({\bf  Case II :} (b) $A(x,t)\ge 0$ and $B(x,t)\ge 0$ ) (b) implies
that
\begin{equation}\rN u(x,t)\le\rP u(x,t).\end{equation}

({\bf  Case III :} (c) $A(x,t)\le 0$ and $B(x,t)\le 0$ ) (c) implies
that $$\rN u(x,t)\ge\rP u(x,t).$$ We note that $-u$ is another
viscosity solution of (8.1). Using $-u$ instead of $u$, we see that
$\rN (-u)(x,t)=\rP u(x,t)$ and $\rP (-u)(x,t)=\rN u(x,t)$. In this
case, the proof can be achieved exactly in the same way as Case II.
Thus we have only to consider Case I and Case II.

Our main goal is to show that there is a universal constant $c>0$
such that $\sup_{Q_r}\rP u\le c\,r^{\ap}$ for any small enough
$r>0$. Since $\rB^{\dd}_r\subset Q_r\subset\rB^{\dd}_{2r}$, it
suffices to show that $\sup_{\rB^{\dd}_r}\rP u\le c\,r^{\ap}$ for
any small enough $r>0$. If we take a rescaled function
$\overline{w}^S_{\vp}(x,t)=\f{1}{c_0} w^S_{\vp}(rx,r^{\sm}t)$ where
$c_0$ is the constant in (8.2), then we may assume that

(i) $|w^S_{\vp}|\le 1$ in $\BR^n_T$ and $\bM^+_2 w^S_{\vp}-\pa_t
w^S_{\vp}\ge-r^\sm/c_0$ in $\rB^{\dd}_1$, for all symmetric sets
$S\subset\BR^n$, and

(ii) for any $(x,t)\in \rB^{\dd}_1$, we have that either
\begin{equation}\f{\ld}{\Ld}\,\rN
u(x,t)-c_1 r^{\sm}(|x|^{\sm}+|t|)^{\f{\bt}{\sm}}\le\rP
u(x,t)\le\f{\Ld}{\ld}\,\rN u(x,t)+c_1
r^{\sm}(|x|^{\sm}+|t|)^{\f{\bt}{\sm}}
\end{equation} or (8.8) holds, for
any small enough $r>0$, where $c_1$ is the constant in (8.7). From
Lemma 3.2, we can also assume that $u$ is $C^{2,\ap_0}$ for some
$\ap_0\in(0,1)$, and so $w^S_{\vp}$, $\rP u$ and $\rN u$ are
continuous.

For our aim, we need only to prove that there are some $r\in(0,1)$
and $\vr\in(0,1)$ such that
\begin{equation}\sup_{\rB^{\dd}_{r^k}}|\rP u|\le(1-\vr)^k=r^{\ap
k}\,\,\,\text{ for $\ap=\f{\ln(1-\vr)}{\ln r}$. }
\end{equation}
We are going to proceed this proof by using mathematical induction.
If $k=0$, then it is trivial by (i). Assume that (8.10) holds in the
$k^{th}$-step ($k\in\BN$). Then we shall show that (8.10) holds also
for the $(k+1)^{th}$-step. By (8.10) and geometric observation, we
have that
\begin{equation}-1\le w^S_{\vp}(x,t)\le\rP u(x,t)\le\f{1}{1-\vr}\,(|x|^{\sm}+|t|)^{\f{\ap}{\sm}}
\end{equation}
for any $(x,t)$ with $(|x|^{\sm}+|t|)^{1/\sm}>r^k$.

We consider the following rescaled functions
\begin{equation*}\begin{split}\widetilde{w}^S_{\vp}(x,t)&:=(1-\vr)^{-k}w^S_{\vp}(r^k
x,r^{k\sm}t),\\
\widetilde{\rP} u(x,t)&:=(1-\vr)^{-k}\rP u(r^k
x,r^{k\sm}t)=\sup_S\widetilde{w}^S_{\vp}(x,t),\\
\widetilde{\rN} u(x,t)&:=(1-\vr)^{-k}\rN u(r^k
x,r^{k\sm}t)=-\inf_S\widetilde{w}^S_{\vp}(x,t).
\end{split}\end{equation*}
Then the function $\widetilde{\rP} u$ satisfies that
\begin{equation*}\widetilde{\rP} u(x,t)\le\begin{cases} 1&\text{ in
$\rB^{\dd}_1$, }\\
\f{1}{1-\vr}\,(|x|^{\sm}+|t|)^{\f{\ap}{\sm}}&\text{ outside
$\rB^{\dd}_1$.} \end{cases}\end{equation*} Choosing $\bt=\ap$ in
(8.9), by (8.8) and (8.9) we have that
\begin{equation}\f{\ld}{\Ld}\,\widetilde{\rN}
u(x,t)-c_1 r^{\sm}\le\widetilde{\rP}
u(x,t)\le\f{\Ld}{\ld}\,\widetilde{\rN} u(x,t)+c_1 r^{\sm}\,\,\text{
in $\rB^{\dd}_1$}\end{equation} and
\begin{equation}\widetilde{\rN} u(x,t)\le\widetilde{\rP} u(x,t)\,\,\text{
in $\rB^{\dd}_1$. }
\end{equation}

Next, we shall show that if $\vr$ and $r$ are chosen so small enough
that $1-\vr=r^{\ap}$ for some $\ap\in(0,1)$, then
$\widetilde{\rP}u\le 1-\vr\,$ in $\rB^{\dd}_r$. This makes it
possible to complete the induction process. For this proof, we
assume that there are some small enough $r$ and $\vr$ such that
$\widetilde{\rP}u\not\le 1-\vr\,$ in $\rB^{\dd}_r$, i.e.
$\widetilde{\rP}u(x_0,t_0)>1-\vr$ for some
$(x_0,t_0)\in\rB^{\dd}_r$. Without loss of generality, we may
suppose that $(x_0,t_0)$ be the point at which the maximum value of
$\widetilde{\rP}u$ is attained in $\rB^{\dd}_r$. Then we see that
\begin{equation}\widetilde{\rP}u(x_0,t_0)=\widetilde{w}^{S_0}_{\vp}(x_0,t_0)>1-\vr
\end{equation} and
\begin{equation}\widetilde{\rP}u(x,t)=\widetilde{w}^{S_0}_{\vp}(x,t)\le\begin{cases}
1&\text{ in
$\rB^{\dd}_1$, }\\
\f{1}{1-\vr}\,(|x|^{\sm}+|t|)^{\f{\ap}{\sm}}&\text{ outside
$\rB^{\dd}_1$,} \end{cases}
\end{equation}
where $S_0$ is the symmetric set given by
$S_0=\{y\in\BR^n:\mu_t(u,x,y)>\mu_0(u,0,y)\}$. Then we note that
\begin{equation}\bM_2^+\widetilde{w}^{S_0}_{\vp}-\pa_t\widetilde{w}^{S_0}_{\vp}
\ge-\f{r^{\sm}}{c_0}\biggl(\f{r^{\sm}}{1-\vr}\biggr)^k>-\f{r^{\sm}}{c_0}
\,\,\,\,\text{ in $\rB^{\dd}_{1/2}$, }\end{equation} because
$\ap<\sm_0<\sm<2$. Since it is easy to check that
$$(1-\widetilde{w}^{S_0}_{\vp})_-\le\bigl(\f{1}{1-\vr}\,(|x|^{\sm}+|t|)^{\f{\ap}{\sm}}-1\bigr)_+:=h(x,t)
\,\,\,\,\text{ in $\BR^n_T$}$$ by (8.15), we derive that
\begin{equation}\bM_2^+(1-\widetilde{w}^{S_0}_{\vp})_-\le\bM_2^+ h
\le c<\iy\,\,\,\,\text{ in $\rB^{\dd}_{1/2}$ }\end{equation} for
some universal constant $c>0$. We also observe that
$$\pa_t(1-\widetilde{w}^{S_0}_{\vp})_-=0\,\,\,\,\text{ in $\rB^{\dd}_{1/2}$, }$$ because
$\rB^{\dd}_1\subset\{(1-\widetilde{w}^{S_0}_{\vp})_-=0\}$ by (8.15).
Let $v^{S_0}_{\vp}=(1-\widetilde{w}^{S_0}_{\vp})_+$. Then we have
that $v^{S_0}_{\vp}(x_0,t_0)=\inf_{\rB^{\dd}_r}v^{S_0}_{\vp}\le\vr$
by (8.14), and moreover by (8.16) and (8.17) we conclude that
\begin{equation*}\begin{split}\bM_2^-
v^{S_0}_{\vp}-\pa_t
v^{S_0}_{\vp}&\le\bM_2^-(1-\widetilde{w}^{S_0}_{\vp})-\pa_t(1-\widetilde{w}^{S_0}_{\vp})\\
&+\bM_2^+(1-\widetilde{w}^{S_0}_{\vp})_- -\pa_t(1-\widetilde{w}^{S_0}_{\vp})_- \\
&\le-(\bM_2^+\widetilde{w}^{S_0}_{\vp}-\pa_t\widetilde{w}^{S_0}_{\vp})\\
&+\bM_2^+(1-\widetilde{w}^{S_0}_{\vp})_-
-\pa_t(1-\widetilde{w}^{S_0}_{\vp})_-\le c\,\,\,\,\text{ in
$\rB^{\dd}_{1/2}$. }
\end{split}\end{equation*}
By Theorem 4.11 \cite{KL3}, there are some universal constants $c>0$
and $\mu>0$ such that
\begin{equation}\bigl|\{v^{S_0}_{\vp}>\ld\vr\}\cap Q_r(x_0,t_0)\bigr|
\le c\,r^{n+\sm}(v^{S_0}_{\vp}(x_0,t_0)+c
r^{\sm})^{\mu}(\ld\vr)^{-\mu}
\end{equation} for any $\ld>0$ and $r\in(0,1/4)$. If we choose $r$
so that $c r^{\sm}<\vr$, then (8.18) becomes
\begin{equation}\bigl|\{v^{S_0}_{\vp}>\ld\vr\}\cap Q_r(x_0,t_0)\bigr|
\le c\,r^{n+\sm}\ld^{-\mu}=c\ld^{-\mu}|Q_r|
\end{equation} for any $\ld>0$. Set $D=\{v^{S_0}_{\vp}\le\ld\vr\}\cap
Q_r(x_0,t_0)$. By (8.19), we have that
\begin{equation}|D|\ge(1-c\ld^{-\mu})|Q_r|\end{equation}
for all large enough $\ld>0$. Since
$v^{S_0}_{\vp}>\ld\vr\,\,\Leftrightarrow\,\,\widetilde{w}^{S_0}_{\vp}<1-\ld\vr$,
we see that $D=\{\widetilde{w}^{S_0}_{\vp}\ge 1-\ld\vr\}\cap
Q_r(x_0,t_0)$. Since $D\subset\rB^{\dd}_1$ and $\widetilde{\rP}u\le
1$ in $D$ by (8.15), we also see that
$\widetilde{\rP}u-\widetilde{w}^{S_0}_{\vp}\le\ld\vr$ in $D$. So we
have the estimate
\begin{equation}\widetilde{\rN}u+\widetilde{w}^{S_0^c}_{\vp}=\widetilde{\rP}u-\widetilde{w}^{S_0}_{\vp}
\le\ld\vr\,\,\text{ in $D$, }
\end{equation}
because
$\widetilde{w}^{S_0}_{\vp}+\widetilde{w}^{S_0^c}_{\vp}=\widetilde{\rP}u-\widetilde{\rN}u$.
For (Case I), it follows from (8.12) and (8.21) that
\begin{equation}\widetilde{w}^{S_0^c}_{\vp}\le-\f{\ld}{\Ld}(1-\ld\vr)+\ld\vr+c_1
r^{\sm}\le-\f{\ld}{2\Ld}\,\,\text{ in $D$, }
\end{equation}
provided that $r$ and $\vr$ are chosen small enough. For (Case II),
by (8.13) and (8.21) we have that
\begin{equation}\widetilde{w}^{S_0^c}_{\vp}\le-(1-\ld\vr)+\ld\vr\le-\f{\ld}{2\Ld}\,\,\text{ in $D$, }
\end{equation}
if $r$ and $\vr$ are chosen small enough. From (8.22), (8.23) and
(8.20), we obtain that
\begin{equation}\bigl|\{\widetilde{w}^{S_0^c}_{\vp}\le-\f{\ld}{2\Ld}\}\cap Q_r(x_0,t_0)\bigr|
\ge(1-c\ld^{-\mu})|Q_r|.
\end{equation}
for any $\ld>0$ and $r\in(0,1/4)$.

For any small $\e>0$, let
$g(x,t)=\bigl(\widetilde{w}^{S_0^c}_{\vp}(r\e
(x-x_0),(r\e)^{\sm}(t-t_0))+\f{\ld}{2\Ld}\bigr)_+$. Then it follows
from (8.24) that
\begin{equation}\bigl|\{g>0\}\cap Q_{\e^{-1}}\bigr|
\le c\ld^{-\mu}|Q_{\e^{-1}}|.
\end{equation}
When $r$ is small enough, by (i) it is also easy to check that
\begin{equation}\bM^+_0 g-\pa_t g\ge-\|u\|_{L^{\iy}_T(L^1_{\om})}\,\,\text{ in $Q_2$.}
\end{equation}
Applying Theorem 5.1 to $g$ with small enough $r\in(0,1/4)$, by
(8.11), (8.15) and (8.25) we obtain that
\begin{equation*}\begin{split}g(x_0,t_0)&\le
C\sup_{s\in(-T,0]}\int_{\BR^n}\f{g(y,s)}{1+|y|^{n+\sm}}\,dy\\
&\le
C\sup_{s\in(-T,0]}\int_{B_{\e^{-1}}}\f{g(y,s)}{1+|y|^{n+\sm}}\,dy
+C\sup_{s\in(-T,0]}\int_{\BR^n\s
B_{\e^{-1}}}\f{g(y,s)}{1+|y|^{n+\sm}}\,dy\\
&\le
C\,\e^{-n-\sm}\ld^{-\mu}+C\,\e^{\ap}\sup_{s\in(-T,0]}\int_{\BR^n\s
B_{\e^{-1}}}\f{|y|^{\ap}+|s|^{\ap/\sm}}{1+|y|^{n+\sm}}\,dy\\
&\le
C\,\e^{-n-\sm}\ld^{-\mu}+\f{C}{\sm-\ap}\,\e^{\sm}+\f{C}{\sm}\,\e^{\sm+\ap}.
\end{split}\end{equation*}
In this estimate, choose $\e$ so small that
$\f{C}{\sm-\ap}\,\e^{\sm}+\f{C}{\sm}\,\e^{\sm+\ap}<\f{\ld}{8\Ld}$,
and then select $\mu$ so large that
$C\,\e^{-n-\sm}\ld^{-\mu}<\f{\ld}{8\Ld}$. Then we have that
$$g(x_0,t_0)\le\f{\ld}{4\Ld}.$$
This implies that
$\widetilde{w}^{S_0^c}_{\vp}(0,0)\le-\f{\ld}{4\Ld}$, which
contradicts to the fact that $\widetilde{w}^{S_0^c}_{\vp}(0,0)=0$.
Hence we conclude that $\widetilde{\rP}u\le 1-\vr\,$ in
$\rB^{\dd}_r$, that is to say, $\rP u\le (1-\vr)^{k+1}\,$ in
$\rB^{\dd}_{r^{k+1}}$. Therefore we complete the proof. \qed

We can also obtain the following corollary in the same manner as
Lemma 8.1.

\begin{cor} If $u\in L^{\iy}_T(L^1_\om)\cap C(\pa_p Q_2)$ be a viscosity solution
satisfying the equation
\begin{equation*}\bI u-\pa_t u=0\,\text{ in
$Q_2$,}\end{equation*} where $\bI$ is defined on $\fL_2(\sm)$ for
$\sm\in(\sm_0,2)$ with $\sm_0\in(1,2)$, then there exists some $\ap\in(0,1)$
such that
$$\sup_{(x,t)\in Q_{1/8}}\f{\rN u(x,t)}{(|x|^{\sm}+|t|)^{\f{\ap}{\sm}}}\lesssim\|u\|_{L^{\iy}_T(L^1_\om)}.$$
\end{cor}

{\bf Proof of Theorem 1.1.} As mentioned above, the case
$\sm\in(0,1]$ could be treated in \cite{KL4}. Thus we have only to
prove our main theorem only for the case $\sm\in(1,2)$.

We note that the fractional Laplacian of order $\sm\in(0,2)$ is
given by
$$-(-\Delta)^{\sm/2}u(x,t)=\int_{\BR^n}\mu_t(u,x,y)\,\f{c_{n,\sm}}{|y|^{n+\sm}}\,dy,$$
where $c_{n,\sm}$ is the constant given below (1.2).
As in (8.4), if $(x,t)\in Q_{1/8}$, then we have that
\begin{equation*}\begin{split}&-(-\Delta)^{\sm/2}u(x,t)+(-\Delta)^{\sm/2}u(0,0)\\
&\qquad=c_{n,\sm}\int_{\BR^n}\bigl[\mu_t(u,x,y)-\mu_0(u,0,y)\bigr]\f{\vp(y)}{|y|^{n+\sm}}\,dy\\
&\qquad\quad+c_{n,\sm}\int_{\BR^n}\bigl[\mu_t(u,x,y)-\mu_0(u,0,y)\bigr]\f{\phi(y)}{|y|^{n+\sm}}\,dy\\
&\qquad=c_{n,\sm}\biggl(\rP u(x,t)-\rN
u(x,t)+\int_{\BR^n}\bigl[\mu_t(u,x,y)-\mu_0(u,0,y)\bigr]\f{\phi(y)}{|y|^{n+\sm}}\,dy\biggr),
\end{split}\end{equation*} where $\vp$ is the radial cut-off function in (8.4) and $\phi=1-\vp$.
Thus it follows from Lemma 8.1, Corollary 8.2 and (8.6) that
$$\sup_{K\in\cK_2}\bigl|C^{\phi}_K(t)\bigr|\lesssim\bigl(|x|^{\sm}+|t|\bigr)^{\f{\ap}{\sm}}\,
\|u\|_{L^{\iy}_T(L^1_{\om})},$$ and thus there is some $\ap\in(0,1)$ such that
\begin{equation}\bigl|(-\Delta)^{\sm/2}u(x,t)-(-\Delta)^{\sm/2}u(0,0)\bigr|\lesssim
(|x|^{\sm}+|t|)^{\f{\ap}{\sm}}\,\|u\|_{L^{\iy}_T(L^1_\om)}
\end{equation} for any $(x,t)\in Q_{1/8}$. Now, by Corollary 3.3,
it is easy to check that
\begin{equation*}\begin{split}\bM^-_2(\btau^t_x u-u)(0,0)&\le\pa_t u(x,t)-\pa_t u(0,0)\\
&=\bI u(x,t)-\bI u(0,0)\le\bM^+_2(\btau^t_x u-u)(0,0).
\end{split}\end{equation*}
Thus, by Lemma 8.1 and Corollary 8.2, we have the estimate
\begin{equation}\begin{split}\bigl|\pa_t u(x,t)-\pa_tu(0,0)\bigr|&\le
|\bM^-_2(\btau^t_x u-u)(0,0)|\vee|\bM^+_2(\btau^t_x u-u)(0,0)|\\
&\le\Ld\bigl(\rP u(x,t)+\rN u(x,t)\bigr)\\
&\lesssim
(|x|^{\sm}+|t|)^{\f{\ap}{\sm}}\,\|u\|_{L^{\iy}_T(L^1_\om)}
\end{split}\end{equation}
for any $(x,t)\in Q_{1/8}$. Hence by a standard translation argument
of (8.27) and (8.28), and the remark (ii) below Theorem 2.1, we
conclude that
$$\|u\|_{C^{\sm+\ap}(Q_{1/2})}\lesssim \,\|u\|_{L^{\iy}_T(L^1_\om)}.$$
Therefore we complete the proof.
\qed

\,\,\noindent{\bf Acknowledgement.} Yong-Cheol Kim was supported by
School of Education, Korea University Grant in 2016.

\end{document}